\documentclass{article}
\usepackage[utf8]{inputenc}
\usepackage{authblk}
\usepackage{setspace}
\usepackage[margin=1.25in]{geometry}
\usepackage{graphicx}
\graphicspath{ {./figures/} }
\usepackage{subcaption}
\usepackage{amsmath}
\usepackage{lineno}
\newtheorem{theorem}{Theorem}
\usepackage{comment}

\usepackage{url}
\usepackage{hyperref}
\usepackage{xcolor}
\hypersetup{
    colorlinks,
    linkcolor={red!50!black},
    citecolor={blue!50!black},
    urlcolor={blue!80!black}
}
\usepackage[natbib=true,backend=biber,sorting=none,style=ieee, citestyle=numeric-comp]{biblatex}

\title{Stochastic Optimal and Time-Optimal Control Studies for Additional Food provided prey-predator Systems involving Holling Type-III Functional  Response}

\author[1]{D Bhanu Prakash}
\author[2]{D K K Vamsi}

\affil[1, 2]{ \ Department of Mathematics and Computer Science, Sri Sathya Sai Institute of Higher Learning, India.}
\affil[1]{First Author. Email: dbhanuprakash@sssihl.edu.in}
\affil[2]{Corresponding author. Email: dkkvamsi@sssihl.edu.in}

\date{}

\begin{document}

\maketitle

\begin{abstract}
This paper consists of a  detailed and novel stochastic optimal control analysis of a coupled non-linear dynamical system. The state equations are modeled as additional food provided prey-predator system with Holling Type-III functional response for predator and intra-specific competition among predators. We firstly discuss the optimal control problem as a Lagrangian problem with a linear quadratic control. Secondly we consider an  optimal control problem in the time-optimal control setting. Stochastic maximum principle is used for establishing the existence of  optimal controls for both these problems.  Numerical simulations are performed based on stochastic forward-backward sweep methods for realizing the theoretical findings.  The results obtained in these optimal control problems are  discussed in the context of biological conservation and pest management.

\end{abstract}


\section{Introduction} \label{Intro}

\qquad Millions of species coexist in this universe where the survival of each species is meticulously woven with respect to other species. The survival of one species(say, predator) is dependent on the existence of other species(say, prey). The very first mathematical models of such interactive species were given by Alfred J. Lotka \cite{lotka1925elements} and Vito Volterra \cite{volterra1927variazioni} in 1925. In the last century, various mathematical models were proposed and their behaviours were studied. \\

One of the  primary component in defining the predator-prey dynamics is the functional response. The functional response is the rate at which each predator captures prey \cite{kot2001elements}. A type III functional response is a sigmoidal response \cite{holling1966functional} that has predators foraging inefficiently at low prey densities. The Holling type III functional response is displayed by many organisms in nature \cite{elkinton2004effects, fernandez2005functional, morozov2010emergence,redpath1999numerical}. Some studies involving additional food-provided prey-predator systems with Holling type III and type IV functional responses can be found in \cite{v1,v2}. Recently, the authors in \cite{ananth2021influence, ananth2022achieving, ananthcmb} have studied controllability  of additional food systems with respect to quality and quantity of additional food  as control variables for type III and type IV systems.
\\

Also, in real world situations, some of the the parameters involved in the model always fluctuate around some average value due to continuous variation in the environment. A large number of researchers introduced a stochastic environmental variation using the Brownian motion into parameters in the deterministic model to construct a stochastic model. The authors in \cite{Type2Stochastic3} proved that a stochastic predator-prey model with a protection zone has a unique stationary distribution which is ergodic. In  \cite{Type3Stochastic} the authors obtained stochastic permenace for a stochastic predator-prey system with Holling-type III functional response. The work \cite{Type2Stochastic4} deals with a  a Holling-type II stochastic predator-prey model with additional food that has an ergodic stationary distribution. The authors in \cite{StochasticType4} studied the deterministic and stochastic dynamics of a modified Leslie-Gower prey-predator system with a simplified Holling-type IV scheme. The work \cite{Type3StochasticAdvanced} deals with the  survival and ergodicity of a stochastic Holling-type III predator-prey model with markovian switching in an impulsive polluted environment. In the recent times, optimal control theory is being applied on various stochastic models in order to achieve optimal control values which minimizes the cost. Also to our knowledge very limited research exists on optimal control theory for  additional food provided stochastic  predator-prey systems involving different functional responses.  \\

Motivated by the above discussions, in this paper, we study two optimal control problems on an additional food provided prey-predator system with Holling type-III functional response for predator. This system is a coupled non-linear dynamical system. The first optimal control problem is a Lagrange problem. This has applications for biological conservation of species where we find the optimal quality and quantity of additional food to be provided to predator in order to maximize the populations of predator and prey \cite{sss1, sss2, sss3}. The second optimal control problem is a special kind of optimal control problem, known as time-optimal control problem, where we determine the optimal additional food to be provided to system in order to reach the final state in minimum time. This has several applications to pest management \cite{ananthcmb, ananth2022achieving, ananth2021influence, srinivasu2011role, varaprasadsirtime}. \\

The rest of the paper is organised as follows. In section \ref{sec2}, we  formulate the stochastic   model for additonal food provided system involving Holling type III response. Section \ref{sec4} deals with the discussion on corresponding linear quadratic optimal control problems with applications to  biological conservation with reference to quality and quantity of additional food as stochastic control variables. Later in   section \ref{sec5} we deal with the time optimal control problems  for  these systems with applications to both biological conservation and  pest management.  Numerical simulations are also done to validate the theoretical findings in sections \ref{sec4}-\ref{sec5}. Finally in section \ref{sec7}, we  do the discussions and conclusions.

\section{Stochastic Model Formulation} \label{sec2}

 In this work we consider the  following deterministic prey-predator model with Holling type-III functional response and additional food for predator given by: 

\begin{equation} \label{det1}
\begin{split}
\frac{\mathrm{d} N(t)}{\mathrm{d} t} & = r N(t) \Bigg(1-\frac{N(t)}{K} \Bigg)-\frac{c N^2(t) P(t)}{a^2+N^2(t)+\alpha \eta  A^2} \\
\frac{\mathrm{d} P(t)}{\mathrm{d} t} & = g \Bigg( \frac{N^2(t) + \eta A^2}{a^2+N^2(t)+\alpha \eta A^2} \Bigg) - m P(t) - d P^2(t)
\end{split}
\end{equation}  \\

The biological descriptions of the various parameters involved in the system (\ref{det1}) are described in Table \ref{param_tab}.  \\

The basic analysis for the model (\ref{det1}) can be done  in  similar lines to analysis done in \cite{srinivasu2018additional}.  \\

To reduce the complexity in the analysis, we now reduce the number of parameters in the model (\ref{det1})  by introducing the transformations $N = a x, P=\frac{a y}{c}$. Then the system (\ref{det1}) gets transformed to:
\begin{equation} \label{det2}
\begin{split}
\frac{\mathrm{d} x(t)}{\mathrm{d} t} & = r x(t) \Bigg( 1 - \frac{x(t)}{\gamma} \Bigg) - \frac{x^2(t) y(t)}{1 + x^2(t) + \alpha \xi } \\
\frac{\mathrm{d} y(t)}{\mathrm{d} t} & = g y(t) \Bigg( \frac{x^2(t) + \xi}{1+x^2(t)+\alpha \xi} \Bigg) - m y(t) - \delta y^2(t)
\end{split}
\end{equation}
where $ \gamma = \frac{K}{a}, \xi = \eta (\frac{A}{a})^2, \delta = \frac{d a}{c}$. \\

As in \cite{Type2Stochastic3, Type3Stochastic}, we now suppose that the intrinsic growth rate of prey and the death rate of predator are mainly affected by environmental noise such that 

$$r \rightarrow r + \sigma_1 dW_1(t), m \rightarrow m+\sigma_2 dW_2(t)$$

where $W_i(t) (i=1,2)$ are the mutually independent standard Brownian motions with $ W_i(0) = 0$ and  $\sigma_1$ and $\sigma_2$ are positive constants and they represent the intensities of the white noise. Hence, the system (\ref{det2}) wiht the environmental noise for parameters $r$ and $m$ and the Holling-type III predator functional response now becomes

\begin{equation} \label{stoc}
\begin{split}
dx & = \left[ r x(t) \Bigg( 1 - \frac{x(t)}{\gamma} \Bigg) - \frac{x^2(t) y(t)}{1 + x^2(t) + \alpha \xi }\right] dt + \sigma_1 x(t) dW_1(t) \\
dy & = \left[ g y(t) \Bigg( \frac{x^2(t) + \xi}{1+x^2(t)+\alpha \xi} \Bigg) - m y(t) - \delta y^2(t) \right] dt + \sigma_2 y(t) dW_2(t)
\end{split}
\end{equation} \\

The biological descriptions of the various parameters involved in the systems  (\ref{det2}) and  (\ref{stoc}) are described in Table \ref{param_tab}.  \\

\begin{table}[bht!]
    \centering
    \caption{Description of variables and parameters present in the systems (\ref{det1}), (\ref{det2}), (\ref{stoc})}

    \begin{tabular}{ccc}
        \hline
        Parameter & Definition & Dimension \\  
        \hline
        T & Time & time\\ 
        N & Prey density & biomass \\
        P & Predator density & biomass \\
        A & Additional food & biomass \\
        r & Prey intrinsic growth rate & time$^{-1}$ \\
        K & Prey carrying capacity & biomass \\
        c & Rate of predation & time$^{-1}$ \\
        a & Half Saturation value of the predators & biomass \\
        g & Conversion efficiancy & time$^{-1}$ \\
        m & death rate of predators in absence of prey & time$^{-1}$ \\
        d & Predator Intra-specific competition & biomass$^{-1}$ time$^{-1}$ \\
        $\alpha$ & quality of additional food & Dimensionless \\
        $\xi$ & quantity of additional food & biomass$^{2}$ \\
        \hline
    \end{tabular}
    \label{param_tab}
\end{table}

\textbf{Existence of global positive solution: } \textit{For any initial value $(x_0,y_0) \in \mathcal{R}_+^2$ there exists a unique solution $(x(t),y(t))$ of system (\ref{stoc}) on $t \geq 0$ and the solution will remain in $\mathcal{R}_+^2$ with probability 1.} \\

The above theorem for existence of solutions of (\ref{stoc}) can be proved in similar lines to the proof in \cite{Type3Stochastic} using the Lyapunov method.

\section{Stochastic Optimal Control Problem} \label{sec4} 

In this section we theoretically establish the existence of optimal control for the system (\ref{stoc}) with both quality and quantity of additional food as stochastic controls respectively.  We also numerically simulate and depict the same with applications to biological conservation.  \\

We now establish the existence of stochastic optimal control using the following stochastic maximum principle \cite{yong1999stochastic}. \\

\newpage

{\large {\bf{Stochastic Maximum Principle:}}}

\begin{theorem}

For a  stochastic controlled system  
\begin{equation} \label{stoc1}
     \begin{cases}
     dX(t) = b(t,X(t),u(t)) dt + \sigma(t,X(t),u(t)) dW(t), t \in [0,T],\\
     X(0)=X_0,
    \end{cases}       
\end{equation}
with the cost functional 
\begin{equation} \label{cost}
J(u(t)) = E \Big[ \int_{0}^{T} f(t,X(t),u(t))  dt + h(X(T)) \Big],    
\end{equation}
with $b(t,X(t),u(t)):[0,T]\times \mathcal{R}^n \times U \rightarrow \mathcal{R}^{n}, \sigma(t,X(t),u(t)):[0,T]\times \mathcal{R}^n \times U \rightarrow \mathcal{R}^{n\times m}$ and for  $T$ denoting the terminal time, $f(t,X(t),u(t))$  the running cost and $h(X(T))$  the terminal cost  with $(\bar{u},\bar{X}) = (\bar{u}, \bar{X}(\bar{u},.))$ as the optimal control and the corresponding optimal state trajectories  \\

then there exist pairs of processes $(p(.),q(.)) \in L_F^2(0,T;\mathcal{R}^n) \times (L_F^2(0,T;\mathcal{R}^n))^m$ satisfying the first order adjoint equations 
\begin{equation} \label{adj}
     \begin{cases}
     dp(t) & = - \Big[ b_X(t,\bar{X}(t),\bar{u}(t))^T p(t) + \sum_{n=1}^{m} \sigma_X^j (t,\bar{X}(t),\bar{u}(t))^T q_j(t) - f_x (t,\bar{X}(t),\bar{u}(t)) \Big] dt + q(t) dW(t),\\
    p(T) & = - h_X(\Bar{X}(T)),    
    \end{cases}       
\end{equation}
and the Variational Inequalitiy (Hamiltonian Maximization Condition)
\begin{equation} \label{max}
    \mathcal{H}(t, \bar{X}(t),\bar{u}(t)) = \max_{u \in U} H(t, \bar{X}(t),u(t),p,q),
\end{equation}
with the Hamiltonian functional given by
\begin{equation} \label{ham}
    H(t,X,u,p,q) = \left\langle p, b \right\rangle + tr[q^{T} \sigma] - f(t,X(t),u(t))
\end{equation} 

\end{theorem}

\subsection{Quality of additional food as control}\label{c1} 

In this section we wish to achieve biological conservation of maximizing prey and predator population for the system (\ref{stoc}) with  quality of additional food as a control variable with minimum supply of the food. \\

To attain this we consider the following objective functional with the state equations (\ref{stoc})

\begin{equation}
 J(u) = E \Bigg[ \int_{0}^{T} \Bigg( -A_1 x(t) - A_2 y(t) + A_3 \frac{\alpha^2(t)}{2} \Bigg) dt \Bigg]  \label{fun}
\end{equation}
where $A_1, A_2, A_3$ are positive constants. \\

Here our goal is to find an optimal control $\alpha^*$ such that $J(\alpha^*) \leq J(\alpha), \forall \alpha(t) \in U$ where $U$ is an admissible control set defined by $U = \{\alpha(t) | 0 \leq \alpha(t) \leq \alpha_{max} \forall t \in (0,t_f]\}$ where $\alpha_{max} \in \mathcal{R}^{+}$. \\

Now to find the optimal control $\alpha^*$ using the stochastic maximum principle we find the similars for the system (\ref{stoc}) with (\ref{stoc1}).
Comparing the stochastic  system (\ref{stoc}) with (\ref{stoc1}), the vectors $X, b, \sigma$ can be seen as

\begin{equation*}
    X(t) = \begin{pmatrix} x(t) \\ y(t) \end{pmatrix}, \  b(t, X(t), u(t)) = \begin{pmatrix} b_1(t) \\ b_2(t) \end{pmatrix} = \begin{pmatrix} r x(t) ( 1 - \frac{x(t)}{\gamma} ) - \frac{x^2(t) y(t)}{1 + x^2(t) + \alpha \xi } \\ g y(t) ( \frac{x^2(t) + \xi}{1+x^2(t)+\alpha \xi} ) - m y(t) - \delta y^2(t) \end{pmatrix}, 
\end{equation*}

\begin{equation*}
    \sigma(t,X(t),u(t)) = \begin{pmatrix} \sigma_1 x(t) & 0 \\ 0 & \sigma_2 y(t) \end{pmatrix},\  dW(t) = \begin{pmatrix} dW_1(t) \\ dW_2(t) \end{pmatrix}
\end{equation*}

We also note that 

\begin{equation} \label{b}
    b_X = \begin{pmatrix} \frac{\partial b_1}{\partial x} & \frac{\partial b_1}{\partial y}\\ \frac{\partial b_2}{\partial x} & \frac{\partial b_2}{\partial y} \end{pmatrix} = \begin{pmatrix} r(1-\frac{2x}{\gamma}) - \frac{2xy(1+\alpha\xi)}{(1+x^2+\alpha \xi)^2} & \frac{-x^2}{1+x^2+\alpha \xi} \\ \frac{2gxy(1+(\alpha-1)\xi)}{(1+x^2+\alpha \xi)^2} & g(\frac{x^2+\xi}{1+x^2+\alpha \xi})-m-2\delta y \end{pmatrix}
\end{equation}

\begin{equation} \label{s}
\begin{split}
\sigma^1 = \begin{pmatrix} \sigma_1 x \\ 0 \end{pmatrix} \Rightarrow \sigma_X^1 = \begin{pmatrix} \frac{\partial }{\partial x} (\sigma_1 x) & \frac{\partial }{\partial y} (\sigma_1 x) \\ \frac{\partial}{\partial x}(0) & \frac{\partial}{\partial y}(0) \end{pmatrix} = \begin{pmatrix} \sigma_1 & 0 \\ 0 & 0 \end{pmatrix}\\
\sigma^2 = \begin{pmatrix} 0 \\ \sigma_2 y \end{pmatrix} \Rightarrow \sigma_X^2 = \begin{pmatrix} \frac{\partial }{\partial x} (0) & \frac{\partial }{\partial y} (0) \\ \frac{\partial}{\partial x}(\sigma_2 y) & \frac{\partial}{\partial y}(\sigma_2 y) \end{pmatrix} = \begin{pmatrix} 0 & 0 \\ 0 & \sigma_2 \end{pmatrix}\\
\end{split} 
\end{equation} \\ 

As the diffusion term $\sigma$ in (\ref{s}) is independent of the control, the solution of the second-order adjoint equations will not be helpful in calculating the optimal control values. \\

Now comparing the cost functional (\ref{fun})  with (\ref{cost}), we see that \\

$f(t,X(t),u(t)) = -A_1 x(t) - A_2 y(t) + A_3 \frac{\alpha^2(t)}{2}$ and $h(X(T)) = 0$. \\

Hence from the stochastic maximum principle,  there exist stochastic processes \\

$(p(.),q(.)) \in L_F^2(0,T;\mathcal{R}^n) \times (L_F^2(0,T;\mathcal{R}^n))^m, \ p(t) = \begin{pmatrix} p_1(t) \\ p_2(t) \end{pmatrix}, \ q(t) = \begin{pmatrix} q_1(t) & q_3(t) \\ q_2(t) & q_4(t) \end{pmatrix} \ni $
\[
\begin{split}
\begin{pmatrix} dp_1(t) \\ dp_2(t) \end{pmatrix} & = - \Bigg[\begin{pmatrix} \frac{\partial b_1}{\partial x} & \frac{\partial b_2}{\partial x}\\ \frac{\partial b_1}{\partial y} & \frac{\partial b_2}{\partial y} \end{pmatrix} \begin{pmatrix} p_1(t) \\ p_2(t) \end{pmatrix} + \sigma_X^1 (t,\bar{X}(t),\bar{u}(t))^T q_1(t) + \sigma_X^2 (t,\bar{X}(t),\bar{u}(t))^T q_2(t) - \begin{pmatrix} \frac{\partial f}{\partial x} \\ \frac{\partial f}{\partial y} \end{pmatrix} \Bigg] dt \\&   + \begin{pmatrix} q_1(t) & q_3(t) \\ q_2(t) & q_4(t) \end{pmatrix} \begin{pmatrix} dW_1(t) \\ dW_2(t) \end{pmatrix}\\
p_1(T) & = 0,  \ p_2(T) = 0.
\end{split}
\]

Substituting $\frac{\partial f}{\partial x} = -A_1, \frac{\partial f}{\partial y} = -A_2 $ and the values of $b_X, \sigma_X^1, \sigma_X^2$ from (\ref{b}), (\ref{s}),  we see that the adjoint equations for the optimal control $\alpha^*$ are given as:

\begin{equation*}
\begin{split}
dp_1(t) &= - \Bigg[ \frac{\partial b_1}{\partial x} p_1(t) + \frac{\partial b_2}{\partial x} p_2(t) + \sigma_1 q_1 + A_1 \Bigg] dt + q_1(t) dW_1(t) + q_3(t) dW_2(t) \\
dp_2(t) &= - \Bigg[ \frac{\partial b_1}{\partial y} p_1(t) + \frac{\partial b_2}{\partial y} p_2(t) + \sigma_2 q_4 + A_2 \Bigg] dt + q_2(t) dW_1(t) + q_4(t) dW_2(t) \\
p_1(T) & = 0,  \ p_2(T) = 0 
\end{split} 
\end{equation*} 

On further simplification we see that 

\begin{equation} \label{adj1}
\begin{split}
dp_1(t) &= - \Bigg[ \Big( r(1-\frac{2x}{\gamma}) - \frac{2xy(1+\alpha\xi)}{(1+x^2+\alpha \xi)^2} \Big) p_1(t) + \frac{2gxy(1+(\alpha-1)\xi)}{(1+x^2+\alpha \xi)^2} p_2(t) + \sigma_1 q_1 + A_1 \Bigg] dt \\& + q_1(t) dW_1(t) + q_3(t) dW_2(t) \\
dp_2(t) &= - \Bigg[ \frac{-x^2}{1+x^2+\alpha \xi} p_1(t) + \Big( g(\frac{x^2+\xi}{1+x^2+\alpha \xi})-m-2\delta y \Big) p_2(t) + \sigma_2 q_4 + A_2 \Bigg] dt \\& + q_2(t) dW_1(t) + q_4(t) dW_2(t) \\
p_1(T) & = 0, \ p_2(T) = 0.
\end{split}
\end{equation}
 
The solutions of the above equations (\ref{adj1}) gives $(p_1(t),p_2(t))$, which are the co-state vectors.  \\

Now  from (\ref{ham}),  we see that the Hamiltonian for the system (\ref{stoc}) is given by: 
\begin{equation*}
\begin{split}
H(t,X,u,p,q) & = \left\langle p, b \right\rangle + tr[q^{T} \sigma] + A_1 x(t) + A_2 y(t) - A_3 \frac{\alpha^2(t)}{2} \\
 & = \begin{pmatrix} p_1(t) & p_2(t) \end{pmatrix} \begin{pmatrix} b_1(t) \\ b_2(t) \end{pmatrix} + tr\Bigg[\begin{pmatrix} q_1(t) & q_2(t) \\ q_3(t) & q_4(t) \end{pmatrix} \begin{pmatrix} \sigma_1 x & 0 \\ 0 & \sigma_2 y \end{pmatrix}\Bigg] + A_1 x(t) + A_2 y(t) - A_3 \frac{\alpha^2(t)}{2} \\
 & = p_1(t) b_1 + p_2(t) b_2 + q_1 \sigma_1 x + q_4 \sigma_2 y + A_1 x(t) + A_2 y(t) - A_3 \frac{\alpha^2(t)}{2} \\
 & = \Big(r x(t) ( 1 - \frac{x(t)}{\gamma} ) - \frac{x^2(t) y(t)}{1 + x^2(t) + \alpha \xi }\Big) p_1(t) + \Big(g y(t) ( \frac{x^2(t) + \xi}{1+x^2(t)+\alpha \xi} ) - m y(t) - \delta y^2(t)\Big) p_2(t) \\&+ q_1 \sigma_1 x + q_4 \sigma_2 y + A_1 x(t) + A_2 y(t) - A_3 \frac{\alpha^2(t)}{2} \\
\end{split}
\end{equation*}

Now from the  Hamiltonian maximization condition (\ref{max}),  we have
\begin{equation*}
\begin{split}
\mathcal{H}(t, \bar{X}(t),\bar{\alpha}(t)) & = \max_{\alpha \in U} H(t, \bar{X}(t),\alpha(t)) \\
\implies \frac{\partial H}{\partial \alpha} & = 0 \\
\end{split}
\end{equation*}

\begin{equation*}
\begin{split}
\implies & \frac{\partial }{\partial \alpha} \Big(r x(t) ( 1 - \frac{x(t)}{\gamma} ) - \frac{x^2(t) y(t)}{1 + x^2(t) + \alpha \xi }\Big) p_1 + \frac{\partial }{\partial \alpha} \Big(g y(t) ( \frac{x^2(t) + \xi}{1+x^2(t)+\alpha \xi} ) - m y(t) - \delta y^2(t)\Big) p_2(t)\\ + & \frac{\partial }{\partial \alpha} (q_1 \sigma_1 x + q_4 \sigma_2 y + A_1 x(t) + A_2 y(t) - A_3 \frac{\alpha^2(t)}{2}) = 0 \\
\implies & \frac{\xi x^2 y p_1}{(1+x^2+\alpha\xi)^2} - \frac{\xi g y p_2 (x^2+\xi)}{(1+x^2+\alpha\xi)^2} - A_3 \alpha = 0 \\ &\ \\
\implies & (A_3 \xi^2) \alpha^3 + (2 A_3 \xi (1+x^2)) \alpha^2 + (A_3 (1+x^2)^2) \alpha + (\xi g y p_2(x^2 + \xi)-\xi x^2 y p_1) = 0.
\end{split}
\end{equation*}

Now from  Descartes' rule of signs, we see that the above cubic equation admits a positive $\alpha$  only if $  g p_2(x^2 + \xi)- x^2 p_1 < 0$.  \\

On solving the above  cubic equation, we see that the optimal quality control is given by 
 
\begin{equation}
    \alpha^* = \frac{(1+x^2-c_0 \xi)^2}{3 c_0 \xi^2},
\end{equation}

\begin{equation*}
\begin{split}
    \text{where, }c_0 &= \bigg(\frac{-10(1+x^2)^3}{\xi^3} + \frac{27y c_1}{2A_3\xi}+\frac{1}{2} \sqrt{\frac{-108 c_1y(1+x^2)^2}{A_3 \xi^4}+\frac{729 y^2c_1^2}{A_3^2\xi^2}}\bigg)^\frac{1}{3}\\
    \text{and }c_1 &= g p_2(x^2+\xi) - x^2 y
\end{split}
\end{equation*} \\

\subsection{Quantity of additional food as control}\label{c2}

In this section we wish to achieve biological conservation of maximizing prey and predator population for the system (\ref{stoc}) with  quantity of additional food as a control variable with minimum supply of the food. \\

To attain this we consider the following objective functional with the state equations (\ref{stoc})

\begin{equation}
 J(u) = E \Bigg[ \int_{0}^{T} \Bigg( -A_1 x(t) - A_2 y(t) + A_3 \frac{\xi^2(t)}{2} \Bigg) dt \Bigg]  \label{fun2}
\end{equation}
where $A_1, A_2, A_3$ are positive constants. \\

Here our goal is to find an optimal control $\xi^*$ such that $J(\xi^*) \leq J(\xi), \forall \xi(t) \in U$ where $U$ is an admissible control set defined by $U = \{\xi(t) | 0 \leq \xi(t) \leq \xi_{max} \forall t \in (0,t_f]\}$ where $\xi_{max} \in \mathcal{R}^{+}$. \\

Now comparing the cost functional (\ref{fun2}) with (\ref{cost}), we see that 

$f(t,X(t),u(t)) = -A_1 x(t) - A_2 y(t) + A_3 \frac{\xi^2(t)}{2}$ and $h(X(T)) = 0$. 

Since $\frac{\partial f}{\partial x} $ and $\frac{\partial f}{\partial y}$ are independent of control parameter, the adjoint equations are same as (\ref{adj1}) in the previous subsection. \\

Hence from the stochastic maximum principle,  there exist stochastic processes \\

$(p(.),q(.)) \in L_F^2(0,T;\mathcal{R}^n) \times (L_F^2(0,T;\mathcal{R}^n))^m, \ p(t) = \begin{pmatrix} p_1(t) \\ p_2(t) \end{pmatrix}, \ q(t) = \begin{pmatrix} q_1(t) & q_3(t) \\ q_2(t) & q_4(t) \end{pmatrix} \ni $

\begin{equation} \label{adj3}
\begin{split}
dp_1(t) &= - \Bigg[ \Big( r(1-\frac{2x}{\gamma}) - \frac{2xy(1+\alpha\xi)}{(1+x^2+\alpha \xi)^2} \Big) p_1(t) + \frac{2gxy(1+(\alpha-1)\xi)}{(1+x^2+\alpha \xi)^2} p_2(t) + \sigma_1 q_1 + A_1 \Bigg] dt \\& + q_1(t) dW_1(t) + q_3(t) dW_2(t) \\
dp_2(t) &= - \Bigg[ \frac{-x^2}{1+x^2+\alpha \xi} p_1(t) + \Big( g(\frac{x^2+\xi}{1+x^2+\alpha \xi})-m-2\delta y \Big) p_2(t) + \sigma_2 q_4 + A_2 \Bigg] dt \\& + q_2(t) dW_1(t) + q_4(t) dW_2(t) \\
p_1(T) & = 0, \ p_2(T) = 0.
\end{split}
\end{equation}
 
The solutions of the above equations (\ref{adj3}) gives $(p_1(t),p_2(t))$, which are the co-state vectors.  \\

Now  from (\ref{ham}),  we see that the Hamiltonian for the system (\ref{stoc}) is given by: 
\begin{equation*}
\begin{split}
H(t,X,u,p,q) & = \left\langle p, b \right\rangle + tr[q^{T} \sigma] + A_1 x(t) + A_2 y(t) - A_3 \frac{\xi^2(t)}{2} \\
 & = \begin{pmatrix} p_1(t) & p_2(t) \end{pmatrix} \begin{pmatrix} b_1(t) \\ b_2(t) \end{pmatrix} + tr\Bigg[\begin{pmatrix} q_1(t) & q_2(t) \\ q_3(t) & q_4(t) \end{pmatrix} \begin{pmatrix} \sigma_1 x & 0 \\ 0 & \sigma_2 y \end{pmatrix}\Bigg] + A_1 x(t) + A_2 y(t) - A_3 \frac{\xi^2(t)}{2} \\
 & = p_1(t) b_1 + p_2(t) b_2 + q_1 \sigma_1 x + q_4 \sigma_2 y + A_1 x(t) + A_2 y(t) - A_3 \frac{\xi^2(t)}{2} \\
 & = \Big(r x(t) ( 1 - \frac{x(t)}{\gamma} ) - \frac{x^2(t) y(t)}{1 + x^2(t) + \alpha \xi }\Big) p_1(t) + \Big(g y(t) ( \frac{x^2(t) + \xi}{1+x^2(t)+\alpha \xi} ) - m y(t) - \delta y^2(t)\Big) p_2(t) \\&+ q_1 \sigma_1 x + q_4 \sigma_2 y + A_1 x(t) + A_2 y(t) - A_3 \frac{\xi^2(t)}{2} \\
\end{split}
\end{equation*}

Now from the  Hamiltonian maximization condition (\ref{max}),  we have
\begin{equation*}
\begin{split}
\mathcal{H}(t, X^*(t),\xi^*(t)) & = \max_{\alpha \in U} H(t, X^*(t),\xi(t)) \\
\implies \frac{\partial H}{\partial \xi} & = 0 \\
\end{split}
\end{equation*}

\begin{equation*}
\begin{split}
\implies & \frac{\partial }{\partial \xi} \Big(r x(t) ( 1 - \frac{x(t)}{\gamma} ) - \frac{x^2(t) y(t)}{1 + x^2(t) + \alpha \xi }\Big) p_1 + \frac{\partial }{\partial \xi} \Big(g y(t) ( \frac{x^2(t) + \xi}{1+x^2(t)+\alpha \xi} ) - m y(t) - \delta y^2(t)\Big) p_2(t)\\ + & \frac{\partial }{\partial \xi} (q_1 \sigma_1 x + q_4 \sigma_2 y + A_1 x(t) + A_2 y(t) - A_3 \frac{\xi^2(t)}{2}) = 0 \\
\implies & \frac{\alpha x^2 y p_1}{(1+x^2+\alpha\xi)^2} + \frac{g y p_2 (1+(1-\alpha)x^2)}{(1+x^2+\alpha\xi)^2} - A_3 \xi = 0 \\
\implies & (A_3 \alpha^2) \xi^3 + (2 A_3 \alpha (1+x^2)) \xi^2 + (A_3 (1+x^2)^2) \xi - (x^2 y p_1 \alpha + g y p_2(1+(1-\alpha)x^2)) = 0
\end{split}
\end{equation*}

Now from  Descartes' rule of signs, we see that the above cubic equation admits a positive $\xi$  only if $x^2 y p_1 \alpha + g y p_2(1+(1-\alpha)x^2) > 0$.  \\

On solving the above  cubic equation, we see that the optimal quantity control is given by 
 
\begin{equation}
    \xi^{*} = \frac{(1+x^2+c_1 \alpha)^2}{3c_1 \alpha^2},
\end{equation}

\begin{equation*}
    \begin{split}
        \text{where, }c_1 &= \Bigg(\frac{-10 (1+x^2)^3}{\alpha^3}-\frac{27 c_2}{2A_3\alpha^2}+\frac{1}{2}\sqrt{\frac{108 (1+x^2)^3 c_2}{A_3 \alpha^5} + \frac{729 c_2^2}{A_3^2 \alpha^4}}\Bigg)^\frac{1}{3},\\
        \text{and, }c_2 &= x^2 y p_1 \alpha + g y p_2(1+(1-\alpha)x^2)
    \end{split}
\end{equation*} \\

\subsection{Numerical Simulations}

In this section, we  numerically illustrate the theoretical findings of the above sections with application to  biological conservation.  \\

Using the Taylor series expansion, the  optimal control problems are simulated and plotted using the Stochastic Forward and Backward Sampling approach. The state equations (\ref{stoc}) and the adjoint equations (\ref{adj1}), (\ref{adj2}) are solved using the forward and backward processes respectively. The forward process is simulated using the Euler-Maruyama scheme \cite{platen2010numerical}. Among the various methods available to discretize the backward process, we chose an implicit scheme with a back propagation of the conditional expectations, which is of order $1/2$ \cite{zhang2004numerical}. These methods are implemented in Python using Sympy, Numpy and Matplotlib packages.  \\

\begin{figure} 
    \begin{subfigure}{0.45\textwidth}
        \includegraphics[width=\textwidth]{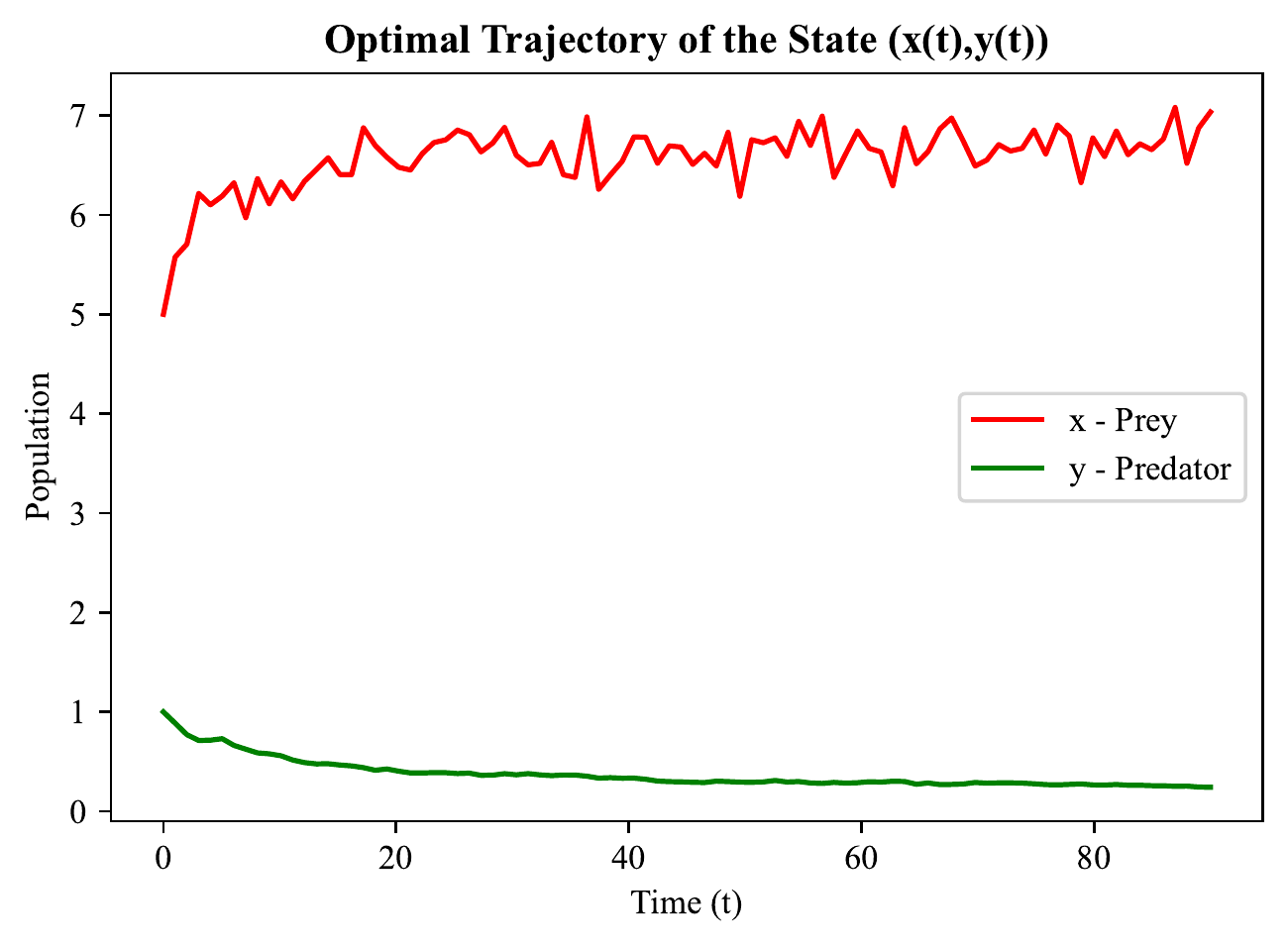}
        \caption{ }
        \label{1a}
    \end{subfigure}
    \begin{subfigure}{0.45\textwidth}
        \includegraphics[width=\textwidth]{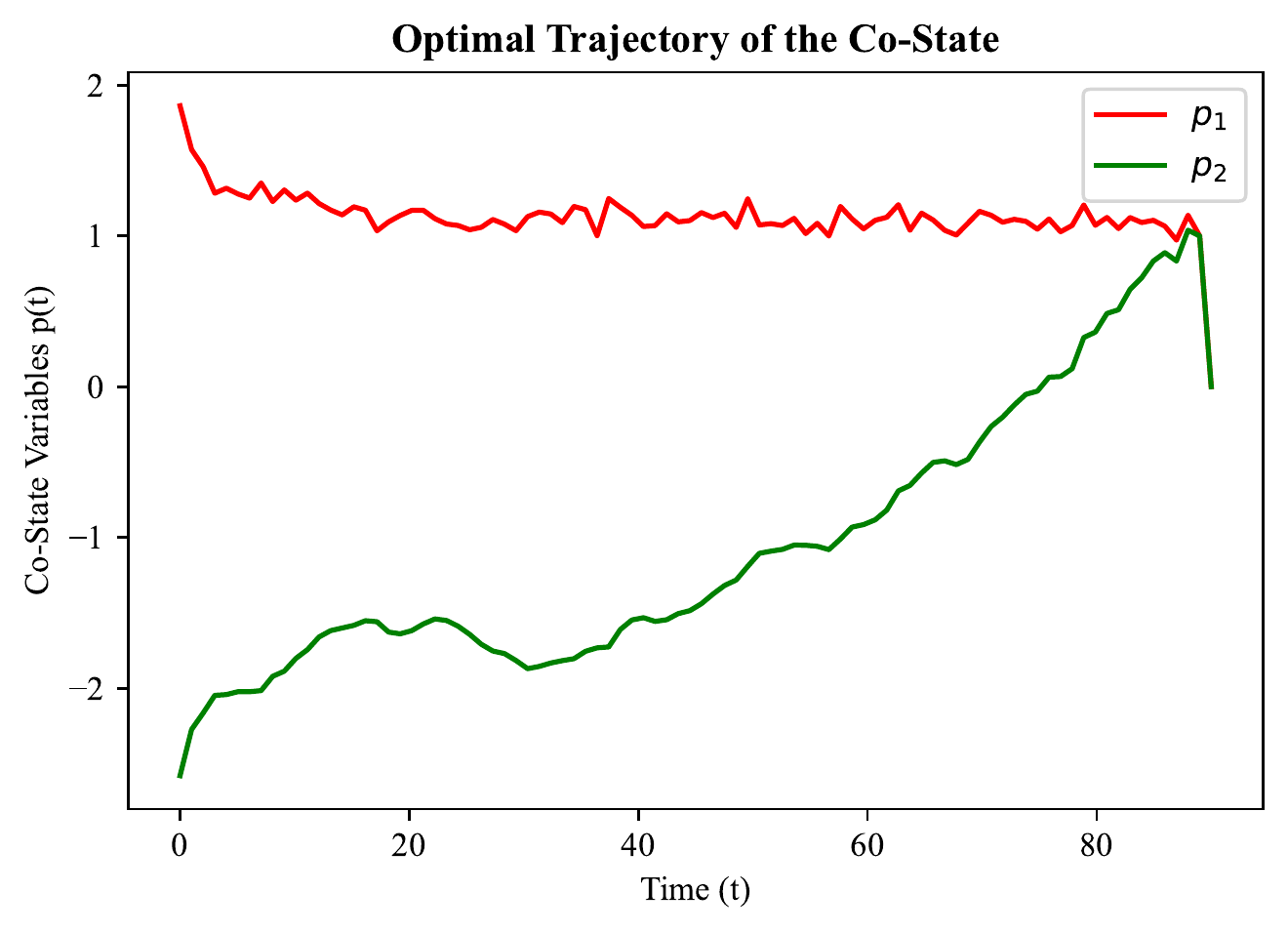}
        \caption{ }
        \label{1b}
    \end{subfigure}
    
    \begin{subfigure}{0.45\textwidth}
        \includegraphics[width=\textwidth]{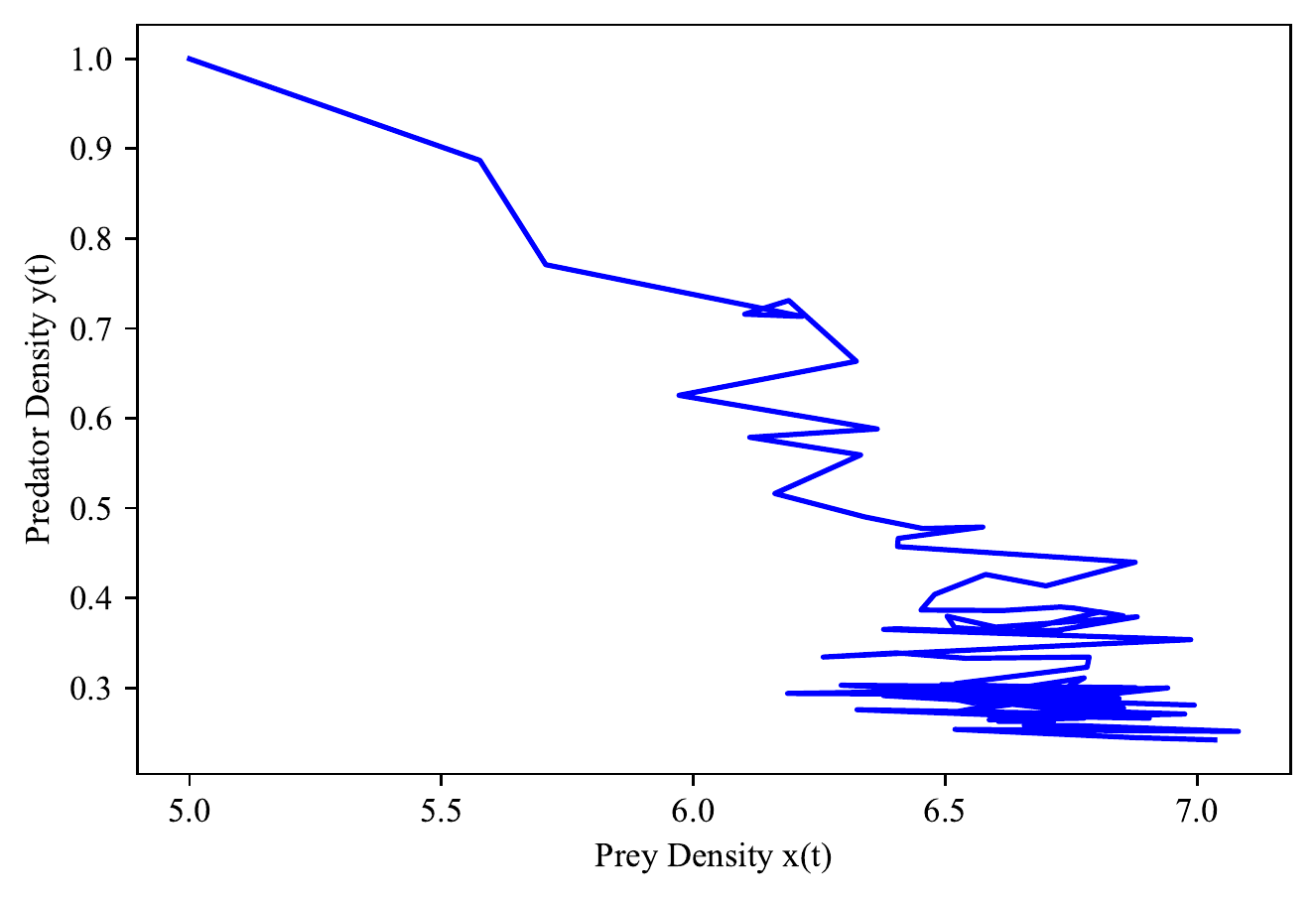}
        \caption{ }
        \label{1c}
    \end{subfigure}
     \begin{subfigure}{0.45\textwidth}
        \includegraphics[width=\textwidth]{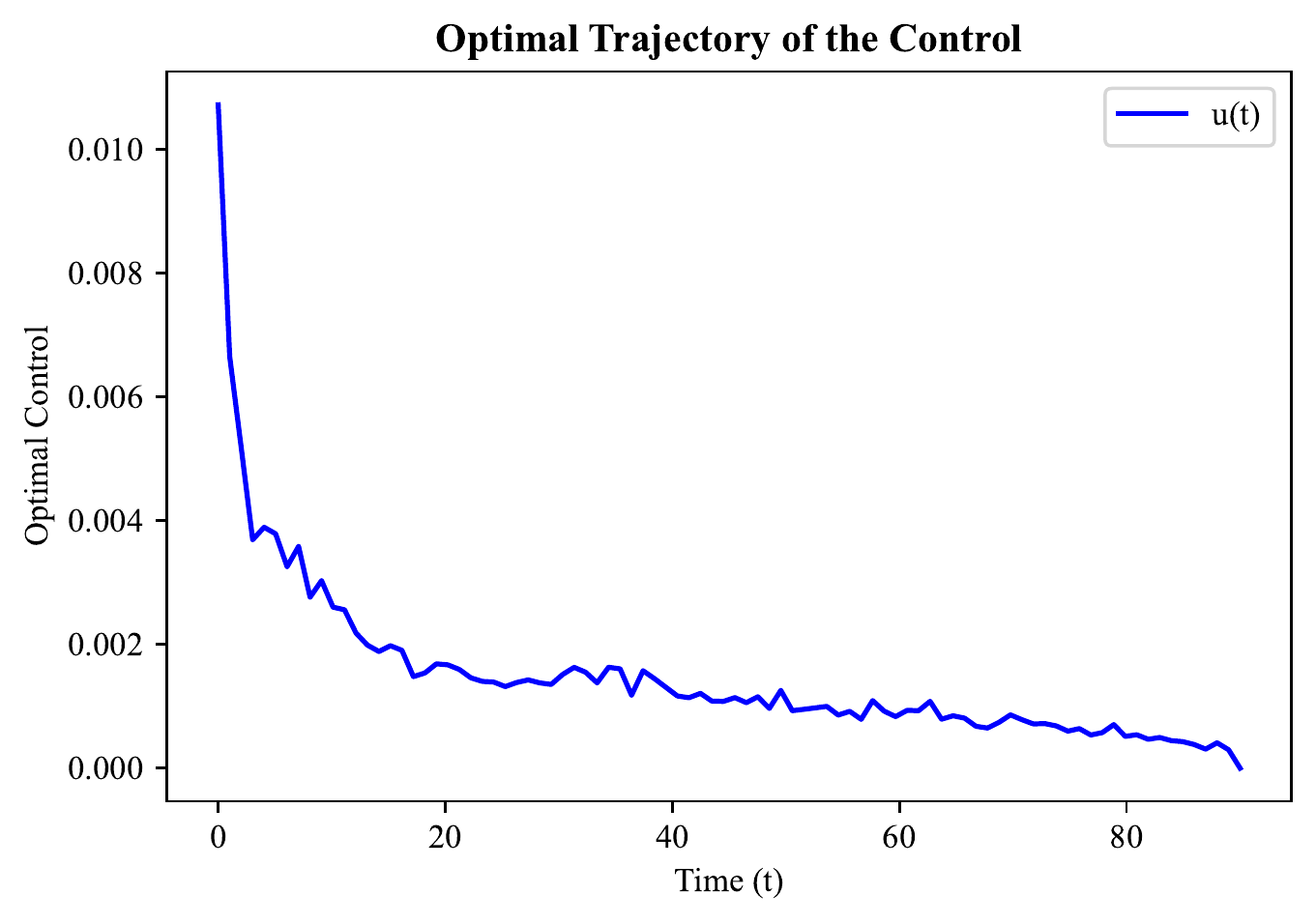}
        \caption{ }
        \label{1d}
    \end{subfigure}

   \caption{This ﬁgure depicts the optimal trajectory of the optimal control problem discussed in section  \ref{c1} with the objective to drive the system from the initial state $(5, 1)$ to the terminal state $(7, 0.25)$. The parameter values for these plots are chosen to be $r=1, \gamma = 7, g = 0.4, m  =0.37, \delta = 0.1, \xi = 0.1, A_1 = A_2 = A_3 = 1$. The terminal value of co-state variables is $(p_1(T), p_2(T)) = (0, 0)$. The intensities of noise are chosen as $\sigma_1 = \sigma_2 = 0.02$. In order to account the randomness, the stochastic control problem is simulated for $5000$ times and the average is plotted in these figures. This example depicts the optimal quality of additional food required to achieve biological conservation, for a fixed quantity of food. The optimal control plot shows that the predators have to be given high quality of additional food initially and then the lower quality of food in order to conserve both the species from extinction. The desired terminal state is reached in $T = 87$ units of time.}
   \label{fig1}
\end{figure}

\begin{figure} 
    \begin{subfigure}{0.45\textwidth}
        \includegraphics[width=\textwidth]{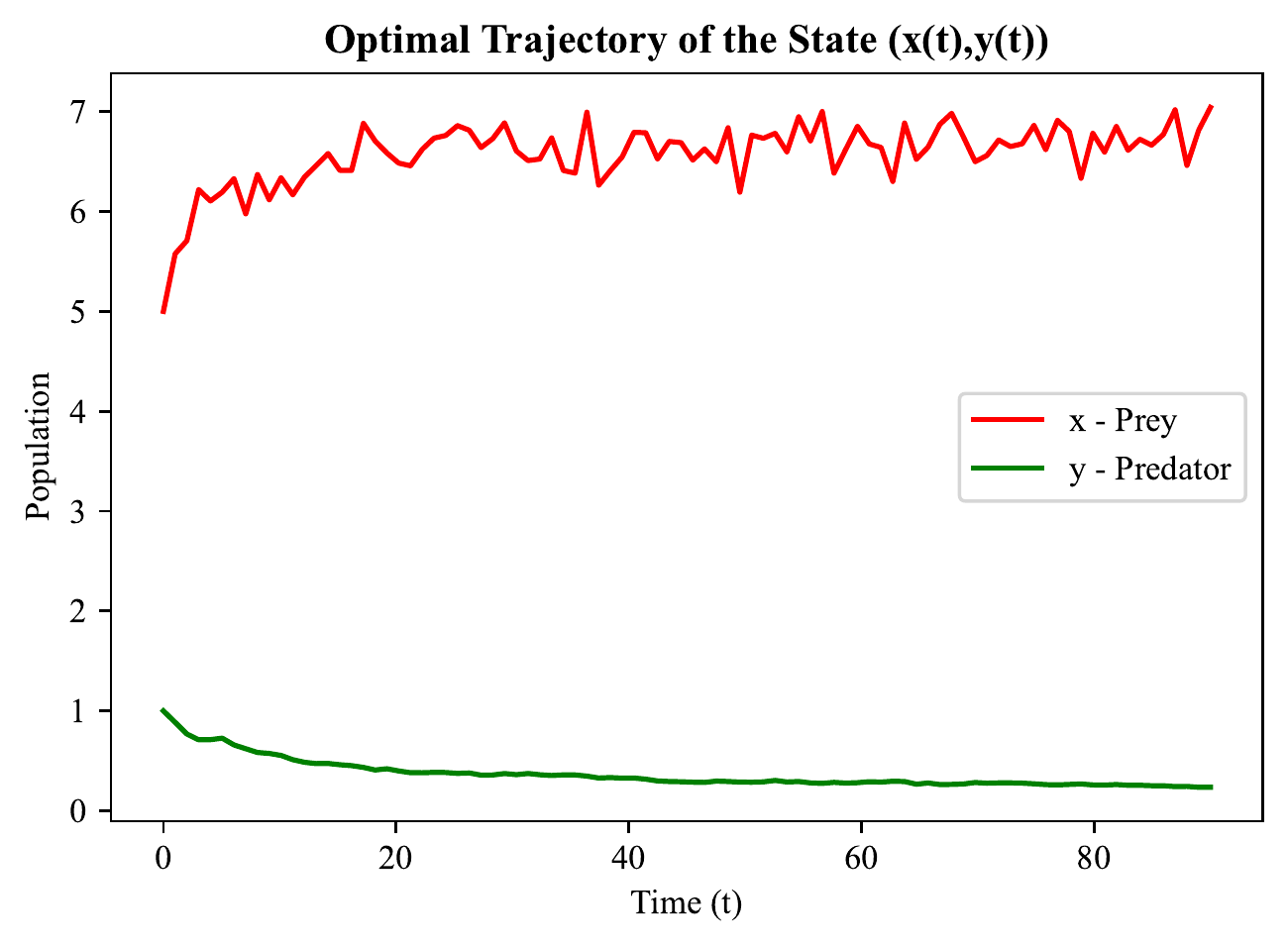}
        \caption{ }
        \label{2a}
    \end{subfigure}
    \begin{subfigure}{0.45\textwidth}
        \includegraphics[width=\textwidth]{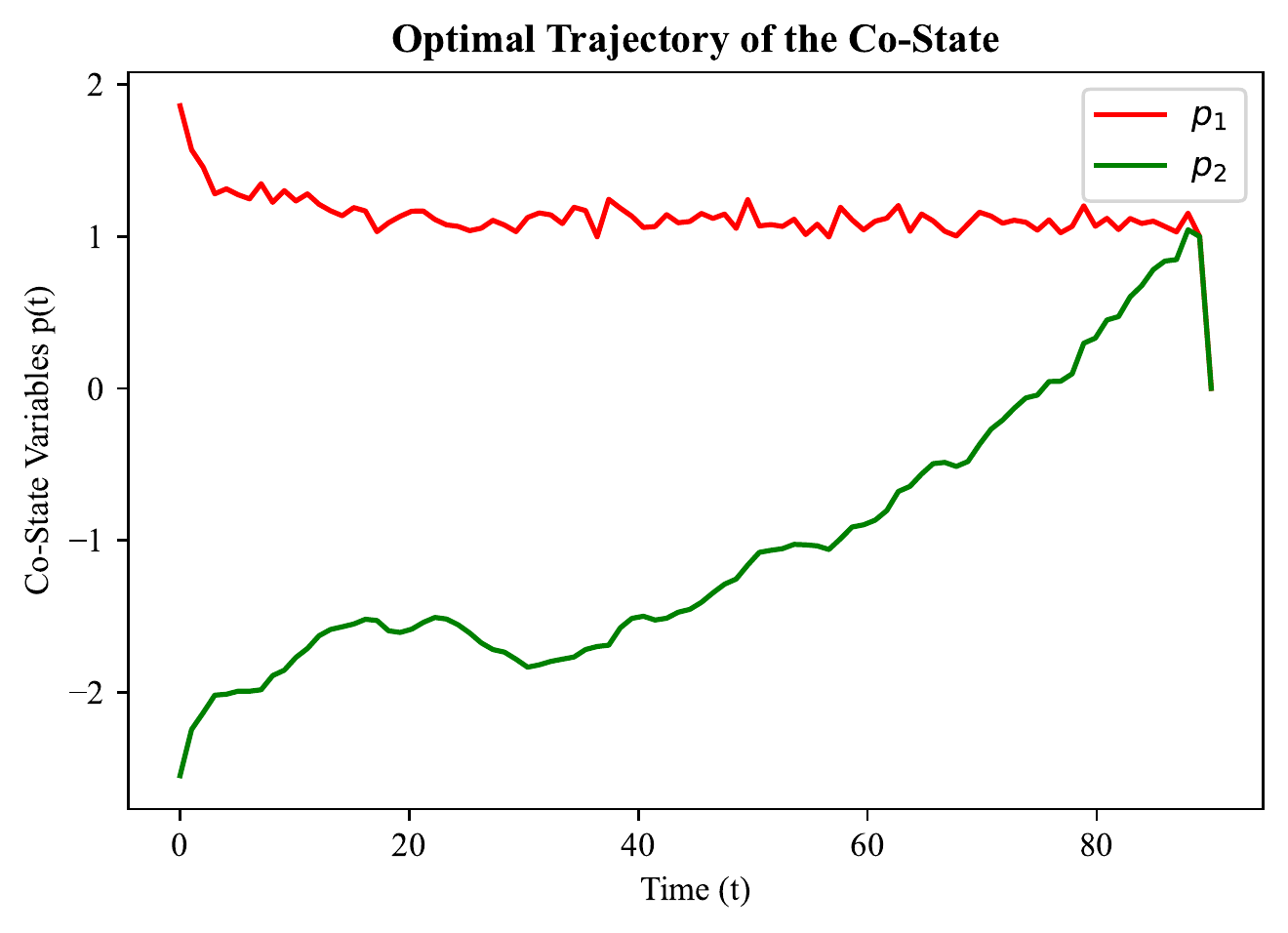}
        \caption{ }
        \label{2b}
    \end{subfigure}
    
    \begin{subfigure}{0.45\textwidth}
        \includegraphics[width=\textwidth]{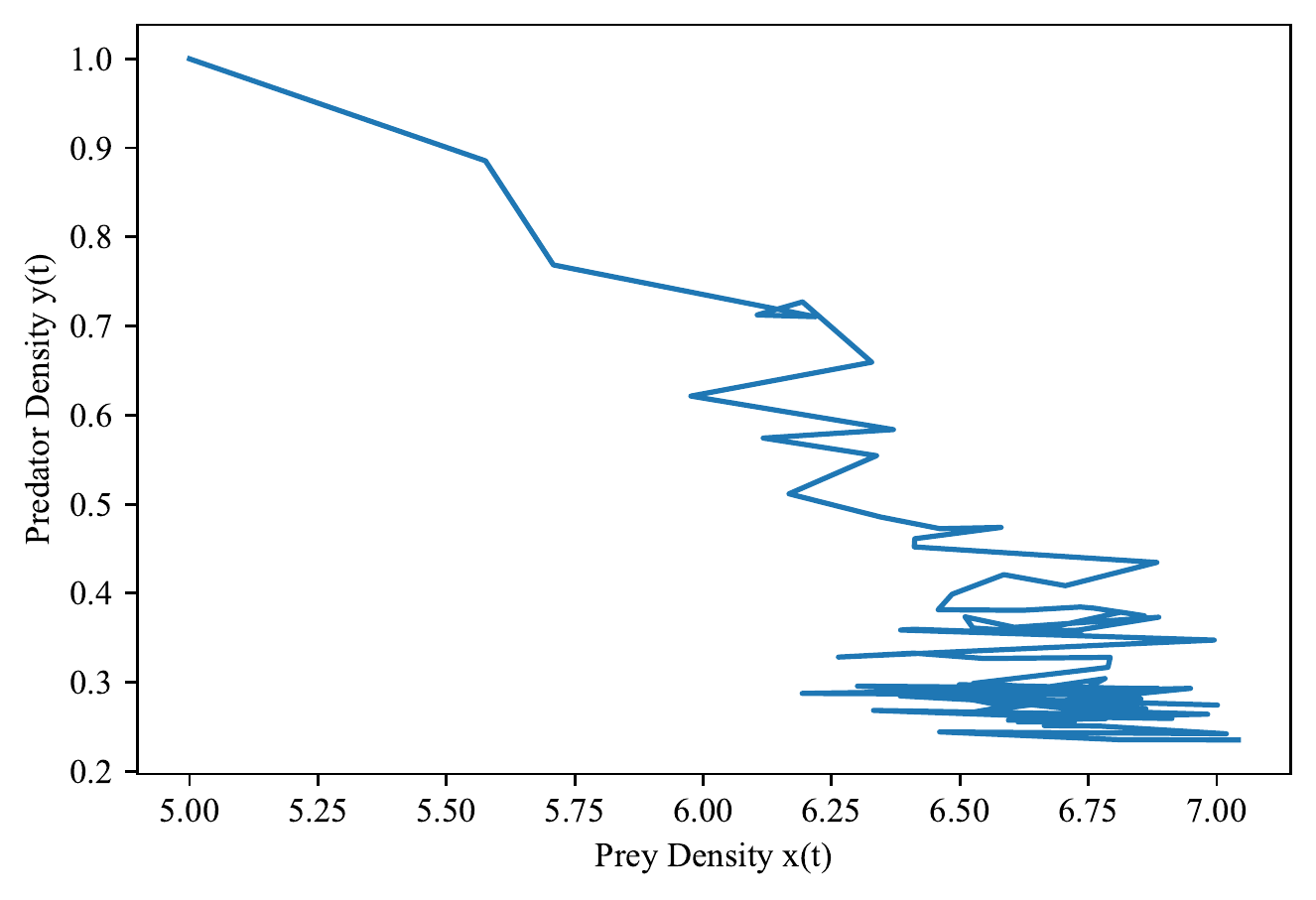}
        \caption{ }
        \label{2c}
    \end{subfigure}
     \begin{subfigure}{0.45\textwidth}
        \includegraphics[width=\textwidth]{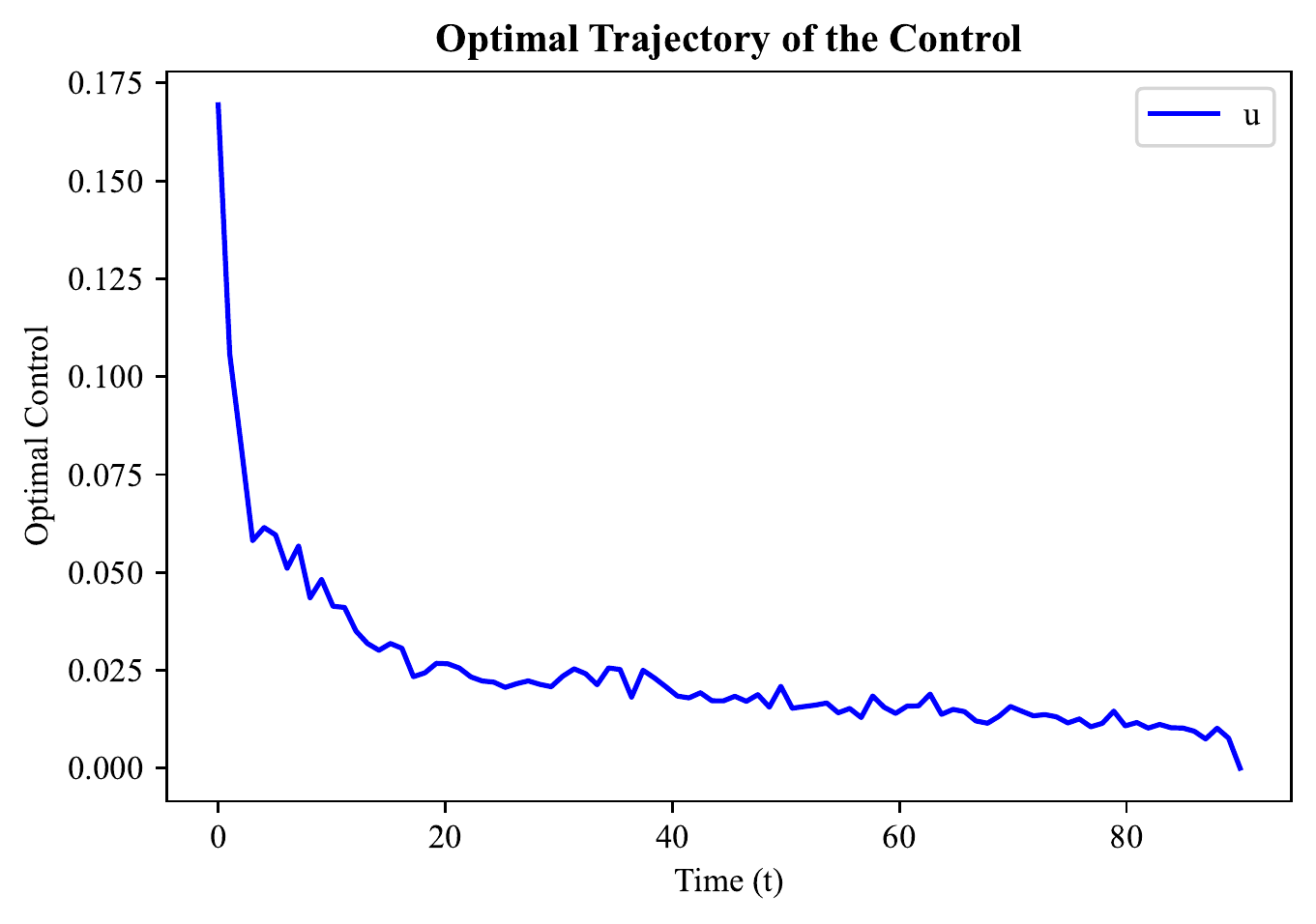}
        \caption{ }
        \label{2d}
    \end{subfigure}

       \caption{This ﬁgure depicts the optimal trajectory of the optimal control problem discussed in section  \ref{c2} with the objective to drive the system from the initial state $(5, 1)$ to the terminal state $(7, 0.25)$. The parameter values for these plots are chosen to be $ r= 1, \gamma = 7, g = 0.4, m = 0.37, \delta = 0.1, \alpha = 2, A_1 = A_2 = A_3 = 1$. The terminal value of co-state variables is $(p_1(T), p_2(T)) = (0, 0)$. The intensities of noise are chosen as $\sigma_1 = \sigma_2 = 0.02$. In order to account the randomness, the stochastic control problem is simulated for $5000$ times and the average is plotted in these figures. This example depicts the optimal quantity of additional food required to achieve biological conservation, for a fixed quality of food. The optimal control plot shows that the predators have to be given high quantity of additional food initially and then the lower quantity of food in order to conserve both the species from extinction. The desired terminal state is reached in $T = 84$ units of time.}
        \label{fig2}
    \end{figure}

The  sub plots present in the below figures \ref{fig1} and \ref{fig2} gives the optimal state trajectories, optimal co-state trajectories, phase diagram and optimal control trajectories respectively. These examples re-iterates the importance of additional food as a control variable in the context of  ecological conservation.

\newpage

\section{Stochastic Time-Optimal Control Problem} \label{sec5}

In this section we wish to achieve pest eradication of nearly prey-elimination stage for the system (\ref{stoc}) with  quality and quantity of additional food as control variables in minimum time. \\

To attain this goal, we consider the following time optimal control probelem with the following objective functional along  with the state equations (\ref{stoc})

\begin{equation}
 J(u) = E \Bigg[ \int_{0}^{T} 1 dt \Bigg].  \label{fun3}
\end{equation}

Here our goal is to find optimal controls $\alpha^*$ and $\xi^*$ such that $J(\alpha^*, \xi^*) \leq J(\alpha, \xi), \forall (\alpha(t), \xi(t)) \in U$ where $U$ is an admissible control set defined by $U = \{(\alpha(t), \xi(t)) | 0 \leq \alpha(t) \leq \alpha_{max}, 0 \leq \xi(t) \leq \xi_{max} \forall t \in (0,t_f]\}$ where $(\alpha_{max},\xi_{max}) \in (\mathcal{R}^{+})^2$\\

Now comparing the cost functional (\ref{fun3})  with (\ref{cost}), we see that \\

$f(t,X(t),u(t)) = 1$ and $h(X(T)) = 0$. \\

Hence from the stochastic maximum principle,  there exist stochastic processes \\

$(p(.),q(.)) \in L_F^2(0,T;\mathcal{R}^n) \times (L_F^2(0,T;\mathcal{R}^n))^m, \ p(t) = \begin{pmatrix} p_1(t) \\ p_2(t) \end{pmatrix}, \ q(t) = \begin{pmatrix} q_1(t) & q_3(t) \\ q_2(t) & q_4(t) \end{pmatrix} \ni $
\[
\begin{split}
\begin{pmatrix} dp_1(t) \\ dp_2(t) \end{pmatrix} & = - \Bigg[\begin{pmatrix} \frac{\partial b_1}{\partial x} & \frac{\partial b_2}{\partial x}\\ \frac{\partial b_1}{\partial y} & \frac{\partial b_2}{\partial y} \end{pmatrix} \begin{pmatrix} p_1(t) \\ p_2(t) \end{pmatrix} + \sigma_X^1 (t,\bar{X}(t),\bar{u}(t))^T q_1(t) + \sigma_X^2 (t,\bar{X}(t),\bar{u}(t))^T q_2(t) - \begin{pmatrix} \frac{\partial f}{\partial x} \\ \frac{\partial f}{\partial y} \end{pmatrix} \Bigg] dt \\&   + \begin{pmatrix} q_1(t) & q_3(t) \\ q_2(t) & q_4(t) \end{pmatrix} \begin{pmatrix} dW_1(t) \\ dW_2(t) \end{pmatrix}\\
p_1(T) & = 0,  \ p_2(T) = 0
\end{split}
\]

Substituting $\frac{\partial f}{\partial x} = 0, \frac{\partial f}{\partial y} = 0 $ and the values of $b_X, \sigma_X^1, \sigma_X^2$ from (\ref{b}), (\ref{s}),  we see that the adjoint equations for the optimal controls $(\alpha^*,\xi^*)$ are given as:

\begin{equation*}
\begin{split}
dp_1(t) &= - \Bigg[ \frac{\partial b_1}{\partial x} p_1(t) + \frac{\partial b_2}{\partial x} p_2(t) + \sigma_1 q_1 \Bigg] dt + q_1(t) dW_1(t) + q_3(t) dW_2(t) \\
dp_2(t) &= - \Bigg[ \frac{\partial b_1}{\partial y} p_1(t) + \frac{\partial b_2}{\partial y} p_2(t) + \sigma_2 q_4 \Bigg] dt + q_2(t) dW_1(t) + q_4(t) dW_2(t) \\
p_1(T) & = 0,  \ p_2(T) = 0 
\end{split}
\end{equation*}

On further simplification, we see that 

\begin{equation} \label{adj2}
\begin{split}
dp_1(t) &= - \Bigg[ \Big(r(1-\frac{2x}{\gamma}) - \frac{2xy(1+\alpha\xi)}{(1+x^2+\alpha \xi)^2} \Big) p_1(t) + \frac{2gxy(1+(\alpha-1)\xi)}{(1+x^2+\alpha \xi)^2} p_2(t) + \sigma_1 q_1 \Bigg] dt \\& + q_1(t) dW_1(t) + q_3(t) dW_2(t) \\
dp_2(t) &= - \Bigg[ \frac{-x^2}{1+x^2+\alpha \xi} p_1(t) + \Big(g(\frac{x^2+\xi}{1+x^2+\alpha \xi})-m-2\delta y \Big)p_2(t) + \sigma_2 q_4 \Bigg] dt \\& + q_2(t) dW_1(t) + q_4(t) dW_2(t) \\
p_1(T) & = 0,  \ p_2(T) = 0
\end{split}
\end{equation}
 
The solutions of the above equations (\ref{adj1}) gives $(p_1(t),p_2(t))$, which are the co-state vectors.  \\

Now  from (\ref{ham}),  we see that the Hamiltonian for the system (\ref{stoc}) is given by: 

\begin{equation*}
\begin{split}
H(t,X,u,p,q) & = \left\langle p, b \right\rangle + tr[q^{T} \sigma]-1 \\
 & = \begin{pmatrix} p_1(t) & p_2(t) \end{pmatrix} \begin{pmatrix} b_1(t) \\ b_2(t) \end{pmatrix} + tr\Bigg[\begin{pmatrix} q_1(t) & q_2(t) \\ q_3(t) & q_4(t) \end{pmatrix} \begin{pmatrix} \sigma_1 x & 0 \\ 0 & \sigma_2 y \end{pmatrix}\Bigg]-1 \\
 & = p_1(t) b_1 + p_2(t) b_2 + q_1 \sigma_1 x + q_4 \sigma_2 y-1\\
 & = \Big(r x(t) ( 1 - \frac{x(t)}{\gamma} ) - \frac{x^2(t) y(t)}{1 + x^2(t) + \alpha \xi }\Big) p_1(t) + \Big(g y(t) ( \frac{x^2(t) + \xi}{1+x^2(t)+\alpha \xi} ) - m y(t) - \delta y^2(t)\Big) p_2(t) \\&+ q_1 \sigma_1 x + q_4 \sigma_2 y - 1 \\
\end{split}
\end{equation*}

Now from the  Hamiltonian maximization condition (\ref{max}),  we have

\begin{equation*}
\begin{split}
\mathcal{H}(t, X^*(t),\alpha^*(t), \xi^*(t)) & = \max_{(\alpha, \xi) \in U} H(t, X^*(t),\alpha(t), \xi(t)) \\
\implies \frac{\partial H}{\partial \alpha} & = 0,  \  \frac{\partial H}{\partial \xi} = 0\\
\end{split}
\end{equation*}
Considering
\begin{equation*}
\begin{split}
\frac{\partial H}{\partial \alpha} & = 0 \\
\implies & \frac{\partial }{\partial \alpha} \Big(r x(t) ( 1 - \frac{x(t)}{\gamma} ) - \frac{x^2(t) y(t)}{1 + x^2(t) + \alpha \xi }\Big) p_1 + \frac{\partial }{\partial \alpha} \Big(g y(t) ( \frac{x^2(t) + \xi}{1+x^2(t)+\alpha \xi} ) - m y(t) - \delta y^2(t)\Big) p_2(t)\\ + & \frac{\partial }{\partial \alpha} (q_1 \sigma_1 x + q_4 \sigma_2 y - 1) = 0 \\
\implies & \frac{\xi x^2 y p_1}{(1+x^2+\alpha\xi)^2} - \frac{\xi g y p_2 (x^2+\xi)}{(1+x^2+\alpha\xi)^2} = 0 \\
\implies & \xi g y p_2(x^2 + \xi)-\xi x^2 y p_1 = 0 \\
\implies & \xi^{*} = \frac{p_1 - g p_2}{g p_2}
\end{split}
\end{equation*}

\noindent
Similarly considering, 
\begin{equation*}
\begin{split}
\frac{\partial H}{\partial \xi} & = 0 \\
\implies & \frac{\partial }{\partial \xi} \Big(r x(t) ( 1 - \frac{x(t)}{\gamma} ) - \frac{x^2(t) y(t)}{1 + x^2(t) + \alpha \xi }\Big) p_1 + \frac{\partial }{\partial \xi} \Big(g y(t) ( \frac{x^2(t) + \xi}{1+x^2(t)+\alpha \xi} ) - m y(t) - \delta y^2(t)\Big) p_2(t)\\ + & \frac{\partial }{\partial \alpha} (q_1 \sigma_1 x + q_4 \sigma_2 y - 1) = 0 \\
\implies & g y p_2(1 + x^2) + \alpha (x^2 y p_1 - g y p_2 x^2) = 0 \\
\implies & \alpha^{*} = \frac{g p_2(1+x^2)}{x^2 (g p_2 - p_1)}
\end{split}
\end{equation*}

\noindent
Hence the optimal quality and quantity variables for the stochastic time optimal control problem are given by  
\begin{center}
$(\alpha^{*}, \xi^*) = \{\frac{g p_2(1+x^2)}{x^2 (g p_2 - p_1)},  \frac{p_1 - g p_2}{g p_2}\}.$
\end{center}

\subsection{Numerical Simulations}

In this section, we  numerically illustrate the theoretical findings of the above time optimal control problem  with application to both biological conservation and pest management. \\

Using the Taylor series expansion, the time optimal control problem is simulated and plotted using the Stochastic Forward and Backward Sampling approach. The state equations (\ref{stoc}) and the adjoint equations (\ref{adj1}), (\ref{adj2}) are solved using the forward and backward processes respectively. The forward process is simulated using the Euler-Maruyama scheme \cite{platen2010numerical}. Among the various methods available to discretize the backward process, we chose an implicit scheme with a back propagation of the conditional expectations, which is of order $1/2$ \cite{zhang2004numerical}. These methods are implemented in Python using Sympy, Numpy and Matplotlib packages.  \\

\begin{figure}
    \begin{subfigure}{0.45\textwidth}
        \includegraphics[width=\textwidth]{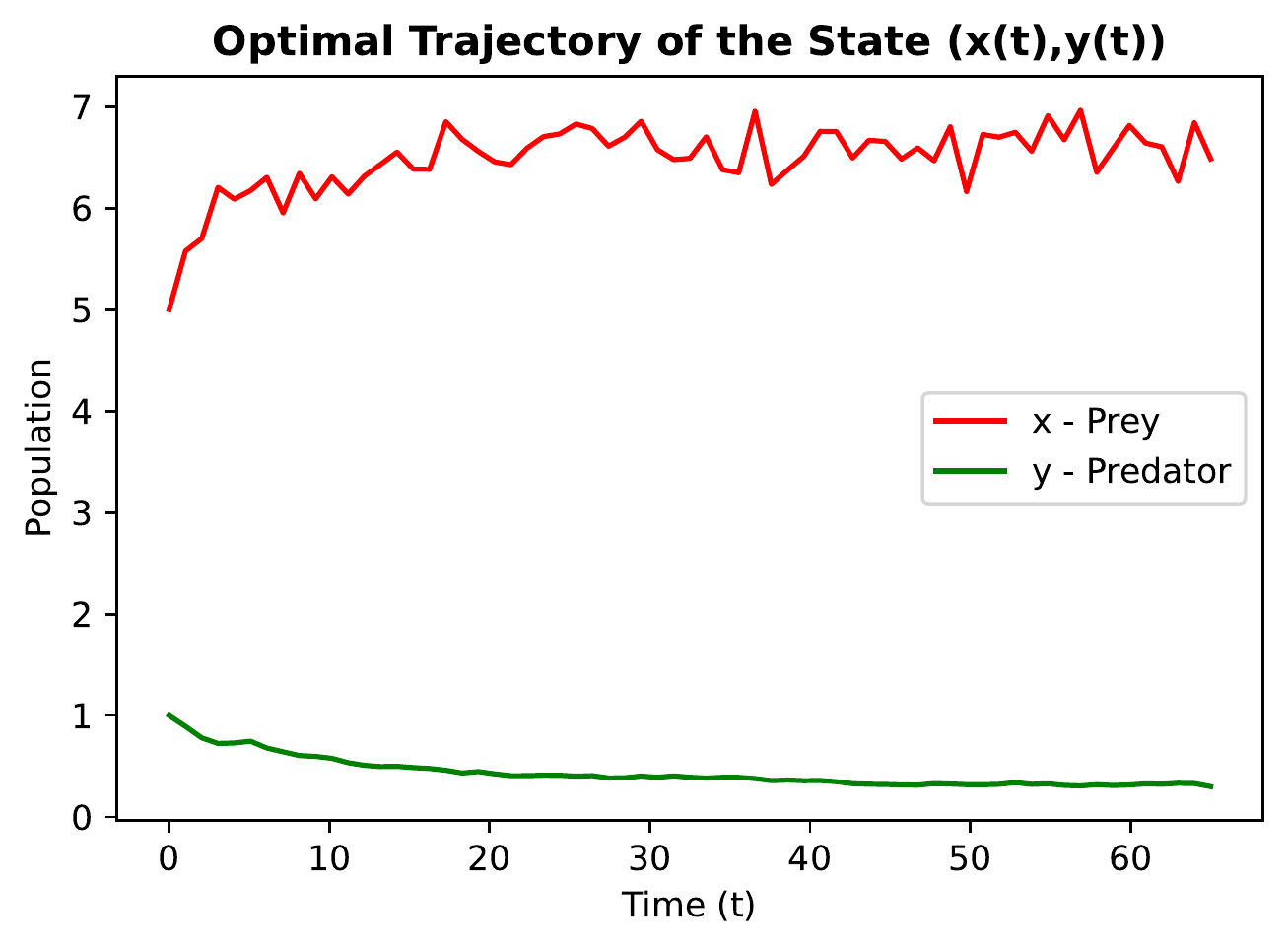}
        \caption{ }
        \label{3a}
    \end{subfigure}
    \begin{subfigure}{0.45\textwidth}
        \includegraphics[width=\textwidth]{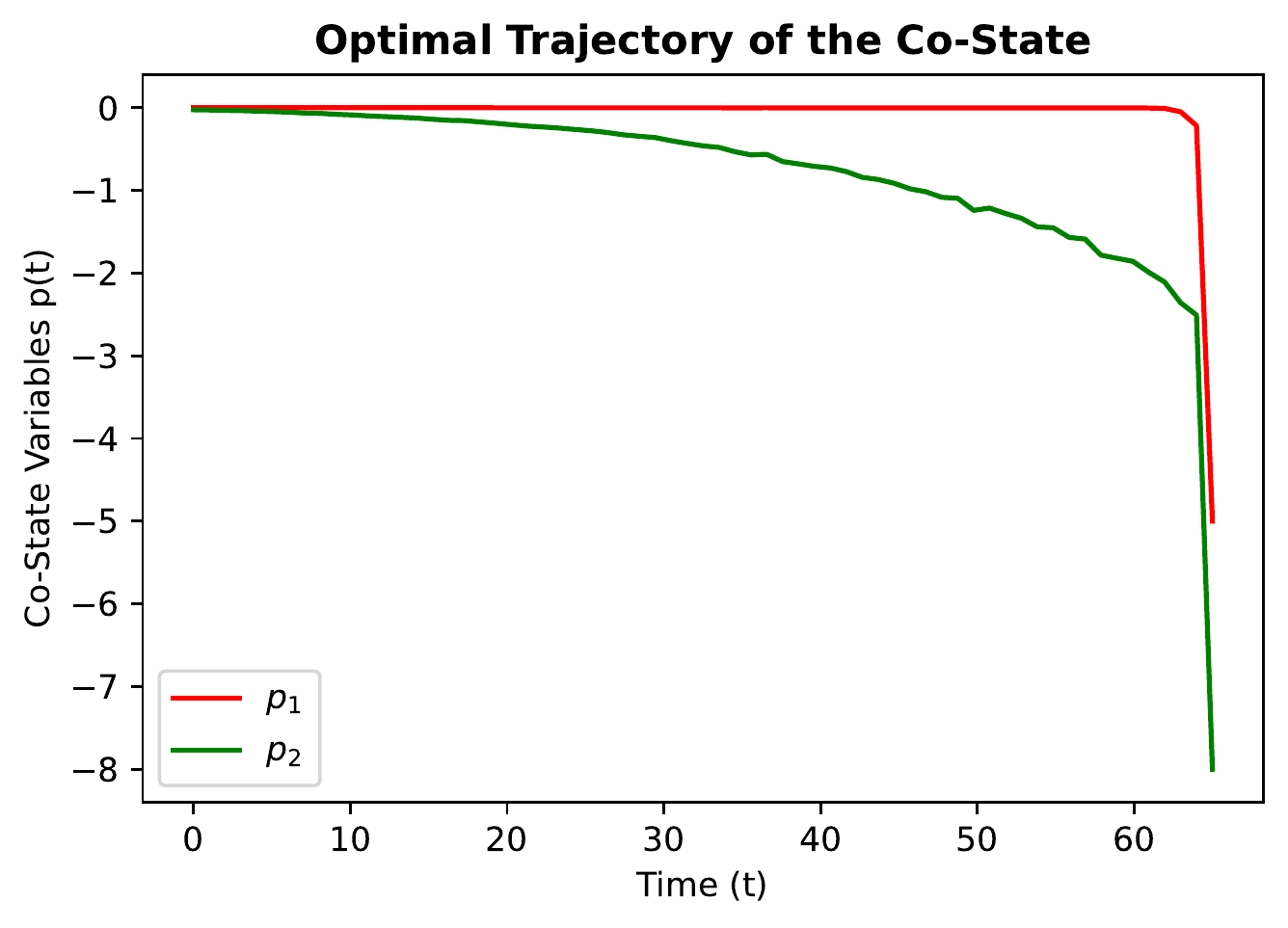}
        \caption{ }
        \label{3b}
    \end{subfigure}
    
    \begin{subfigure}{0.45\textwidth}
        \includegraphics[width=\textwidth]{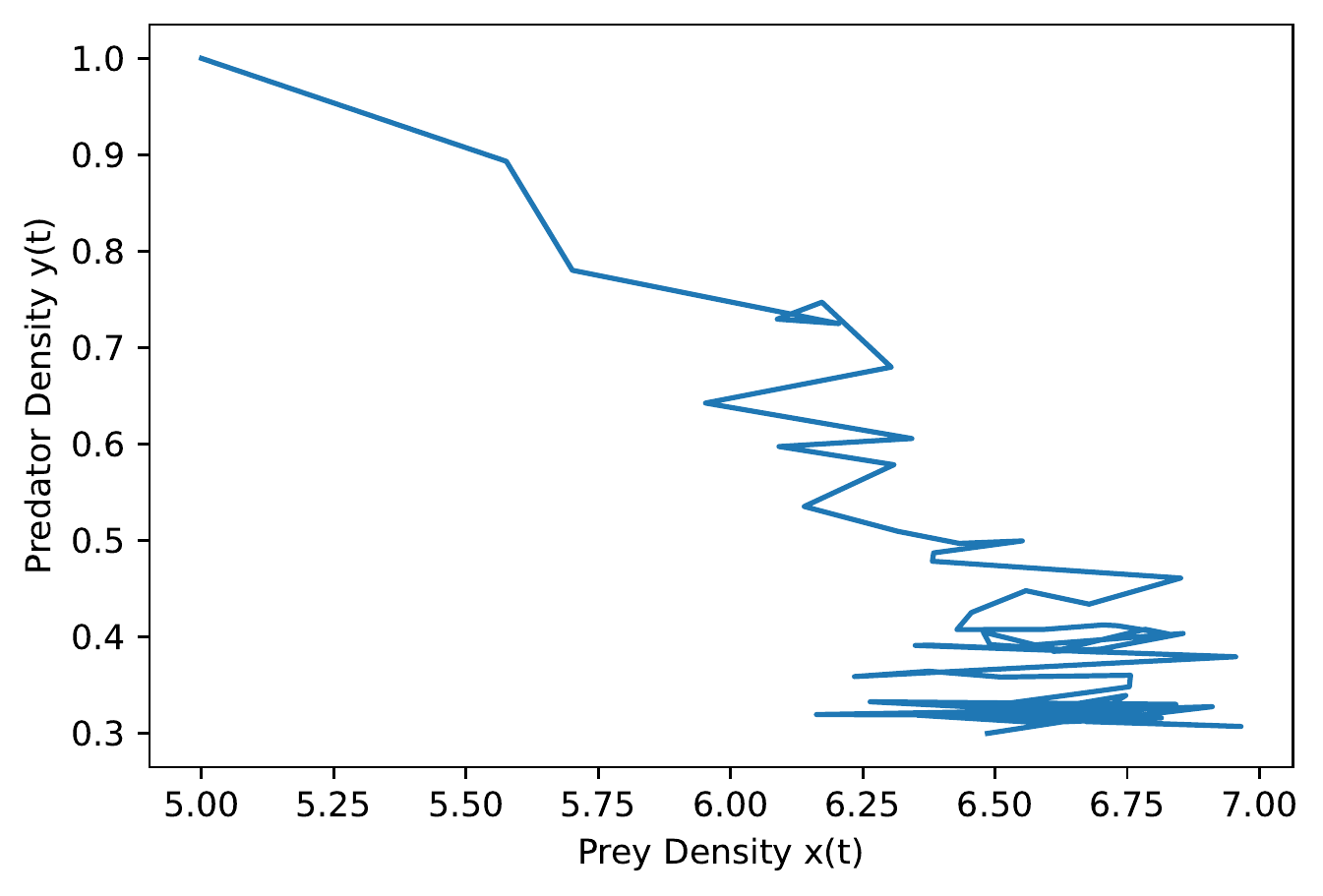}
        \caption{ }
        \label{3c}
    \end{subfigure}
     \begin{subfigure}{0.45\textwidth}
        \includegraphics[width=\textwidth]{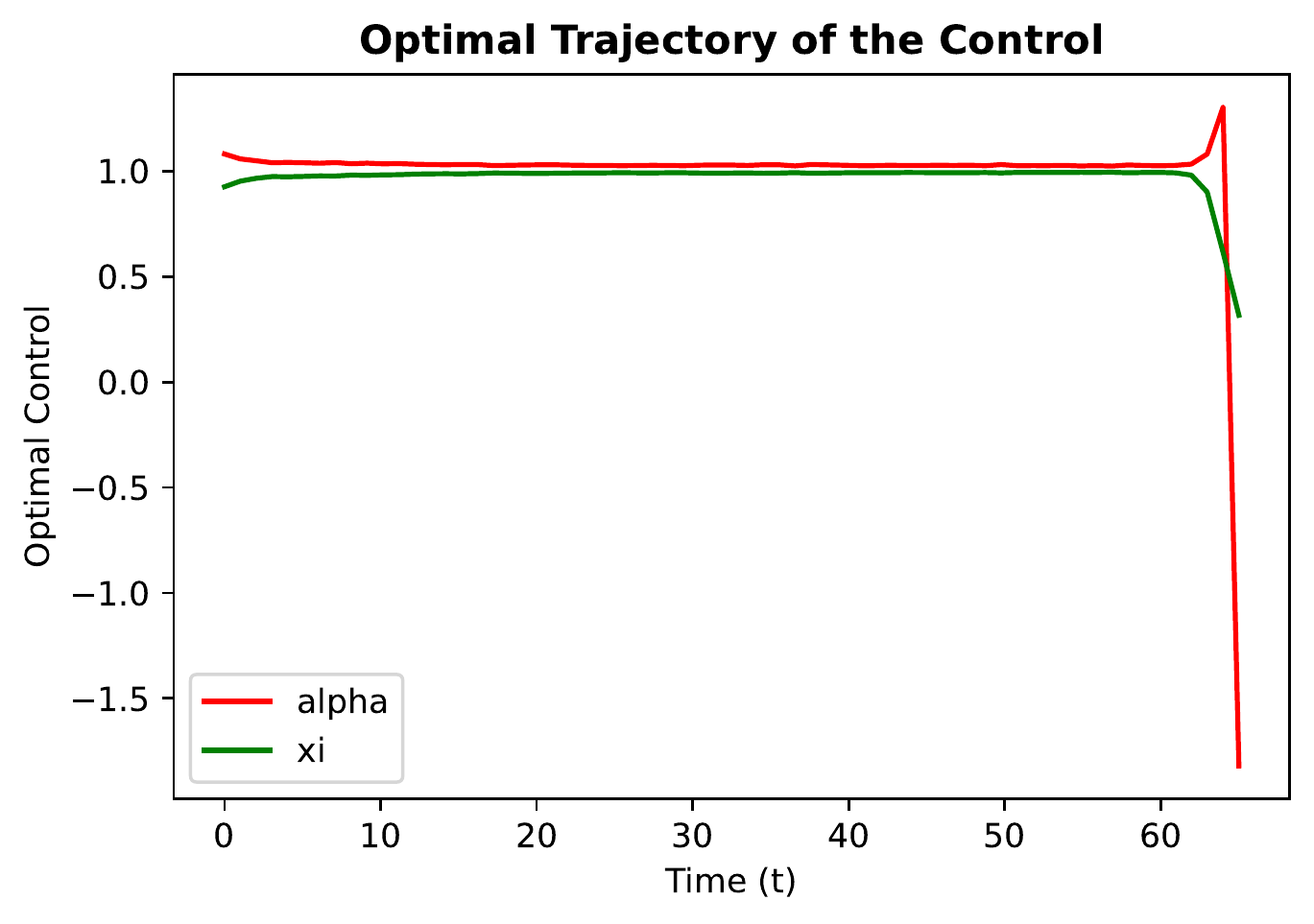}
        \caption{ }
        \label{3d}
    \end{subfigure}

    \caption{This ﬁgure depicts the optimal trajectory of the time-optimal control problem discussed in section \ref{sec5} with the objective to drive the system from the initial state $(5, 1)$ to the terminal state $(6.75, 0.3)$ in minimum time. The parameter values for these plots are chosen to be  $ r = 1, \gamma = 7, g = 0.4, m = 0.37, \delta = 0.1, \alpha = 2, \xi = 0.5$. The initial value of co-state variables is $(p_1(0), p_2(0)) = (0, 0)$. The intensities of noise are chosen as $\sigma_1 = \sigma_2 = 0.02$. In order to account the randomness, the stochastic time-optimal control problem is simulated for $5000$ times and the average is plotted in these figures. This example depicts the optimal quality and quantity of additional food required to achieve biological conservation. The optimal control plot shows that the predators have to be given constant quantity of additional food ($1$ unit) and the constant quantity of additional food ($1$ unit) in order to conserve both the species from extinction. The desired terminal state is reached in $T = 63$ units of time.}
     \label{fig3}
\end{figure}

\begin{figure} 
    \begin{subfigure}{0.45\textwidth}
        \includegraphics[width=\textwidth]{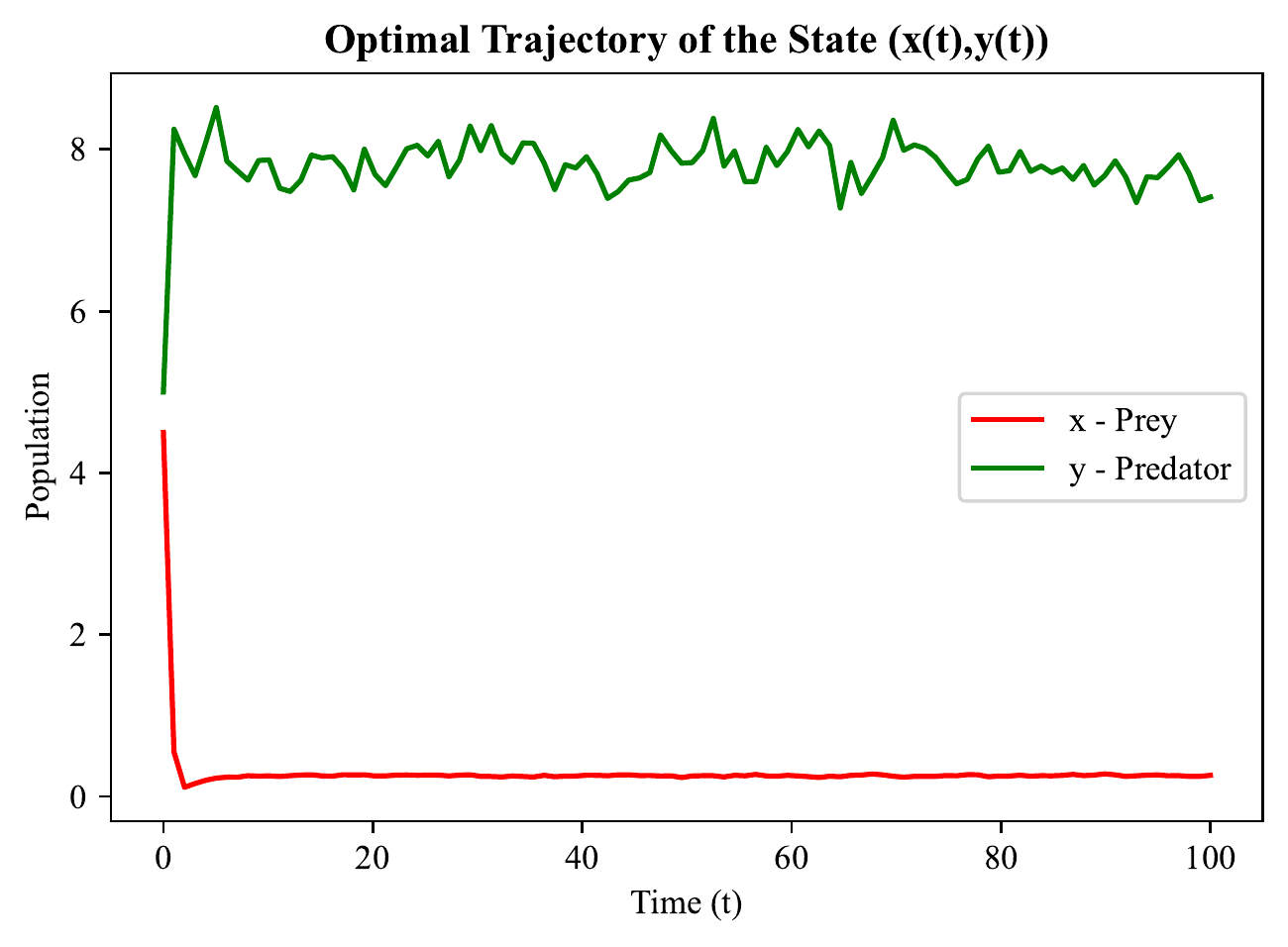}
        \caption{ }
        \label{4a}
    \end{subfigure}
    \begin{subfigure}{0.45\textwidth}
        \includegraphics[width=\textwidth]{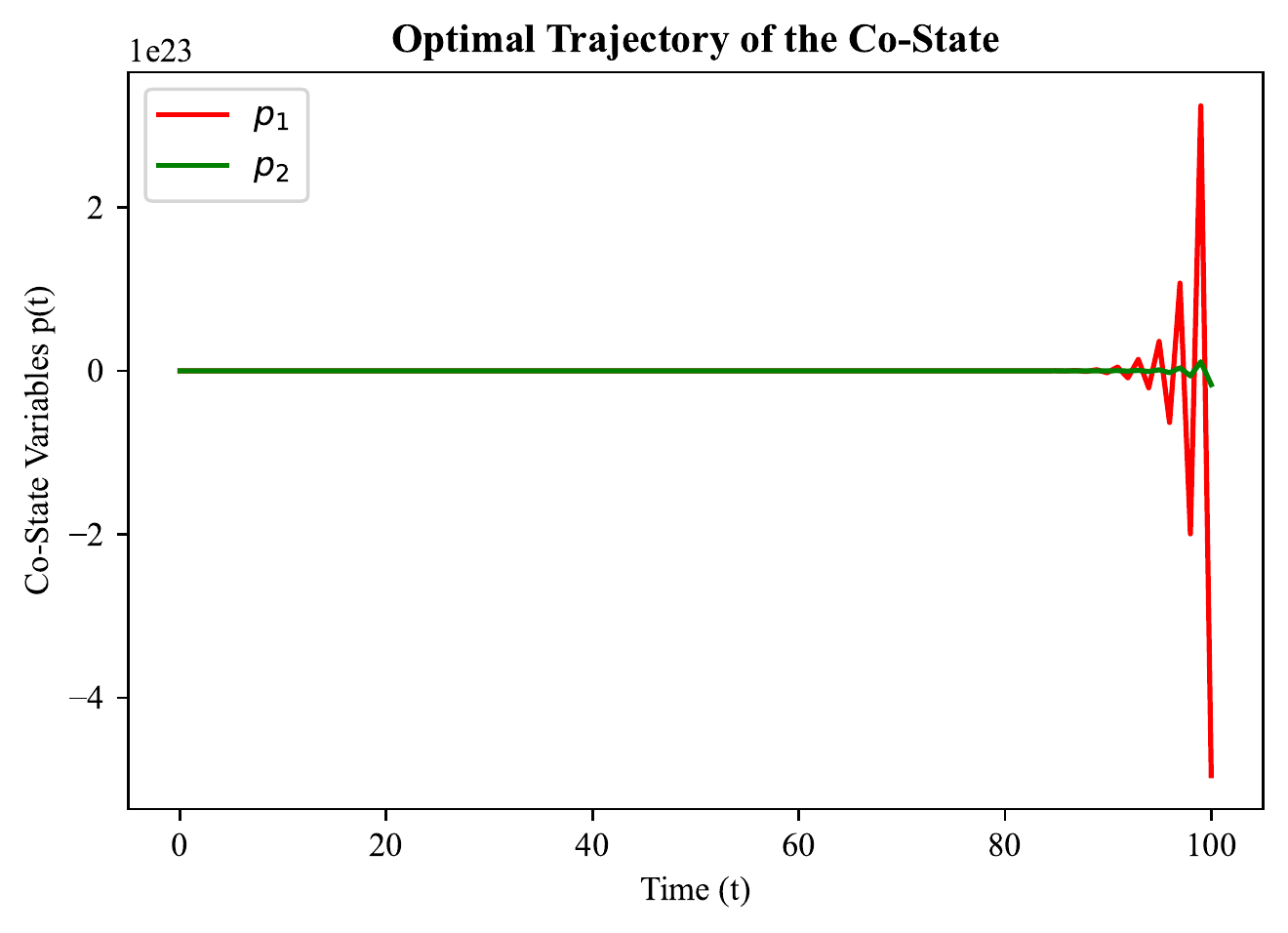}
        \caption{ }
        \label{4b}
    \end{subfigure}
    
    \begin{subfigure}{0.45\textwidth}
        \includegraphics[width=\textwidth]{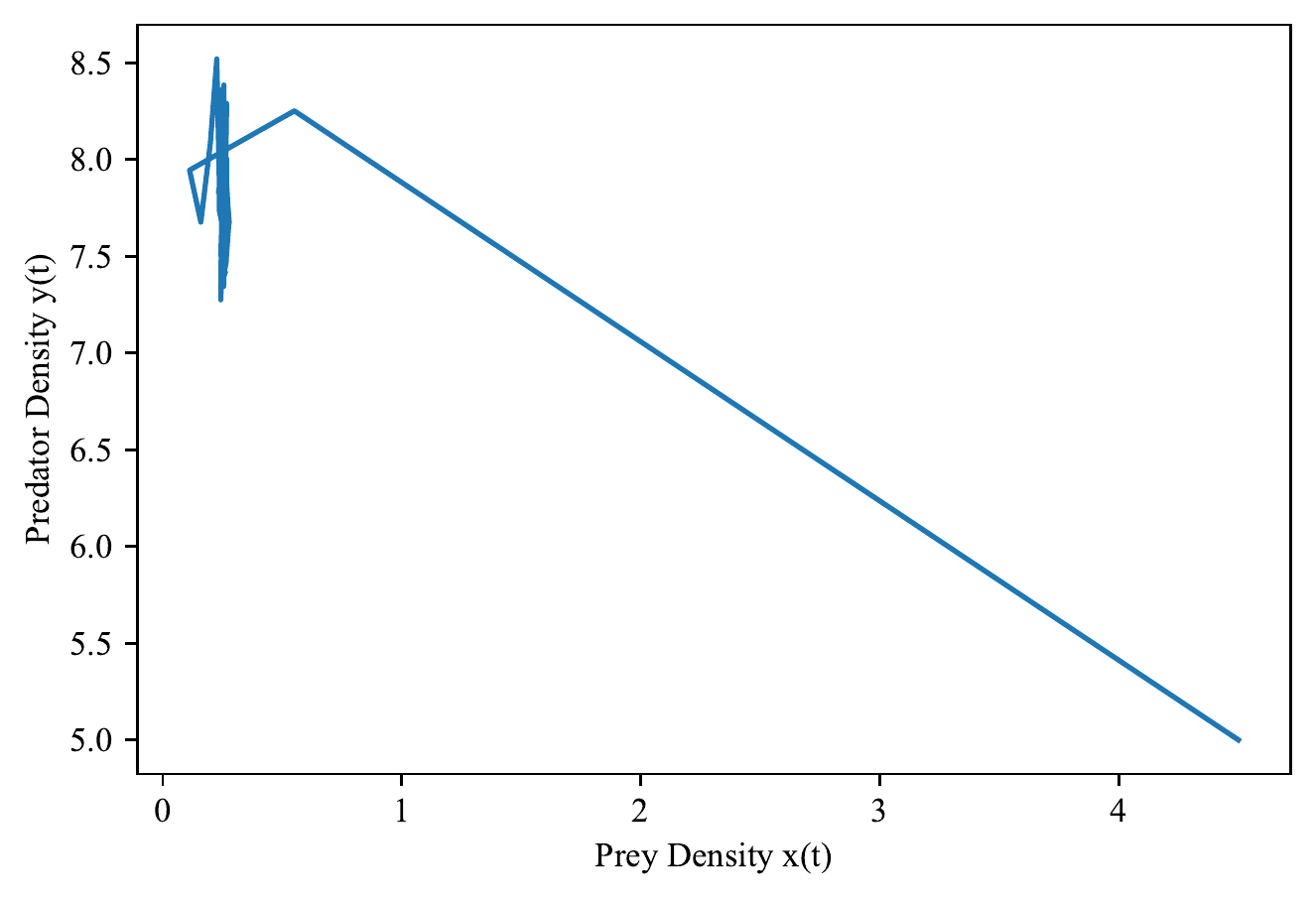}
        \caption{ }
        \label{4c}
    \end{subfigure}
     \begin{subfigure}{0.45\textwidth}
        \includegraphics[width=\textwidth]{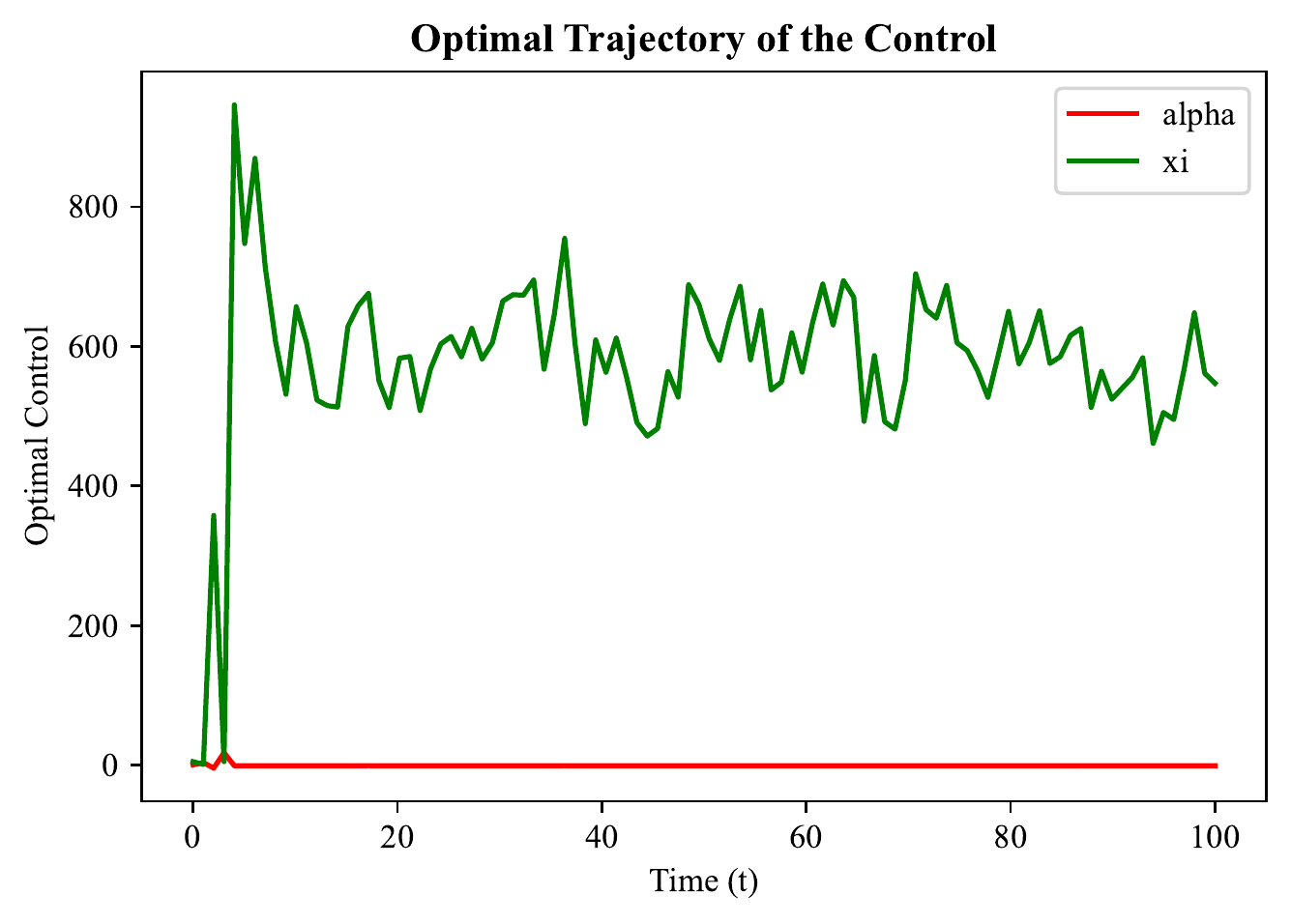}
        \caption{ }
        \label{4d}
    \end{subfigure}

    \caption{This ﬁgure depicts the optimal trajectory of the time-optimal control problem discussed in section \ref{sec5} with the objective to drive the system from the initial state $(4.5, 5)$ to the terminal state $(0.26, 7.41)$ in minimum time. The parameter values for these plots are chosen to be $ r = 0.7, \gamma = 5, g = 1.2, m = 0.19, \delta = 0.05, \alpha = 1.4, \xi = 1.4$. The initial value of co-state variables is $(p_1(0), p_2(0)) = (9, -6)$. The intensities of noise are chosen as $\sigma_1 = \sigma_2 = 0.02$. In order to account the randomness, the stochastic time-optimal control problem is simulated for $5000$ times and the average is plotted in these figures. This example depicts the optimal quality and optimal quantity of additional food required to achieve pest management. The optimal control plot shows that the predators have to be given constant quality of additional food ($0.67$ units) and a high  quantity of additional food   throughout for reducing the pest population effectively. This example suggest  that even with   low quality of additional food and high quantity, the prey (pest) can be reduced to a threshold level at which they no longer cause signiﬁcant damage to the ecosystem. The desired terminal state is reached in $T = 100$ units of time.}
    \label{fig4}
\end{figure}

The sub plots present in the below figures \ref{fig3} and \ref{fig4} gives the optimal state trajectories, optimal co-state trajectories, phase diagram and optimal control trajectories respectively. These examples re-iterates the importance of additional food as control variables in the context of both  ecological conservation and pest management.

\newpage
\section{Discussions and Conclusions} \label{sec7}

The provision of additional food has proven to be very effective in conserving endangered species \cite{harwood2004prey, putman2004supplementary, redpath2001does} as well as controlling invasive or harmful species \cite{sabelis2006does, wade2008conservation, winkler2005plutella}. A significant amount of theoretical work was also done on the impact of additional food in various ecosystems \cite{varaprasadsiraddi,varaprasadsirtime}.  Both theoretically and experimentally, the technique of providing additional food for bio-control seemed to be very effective. Now a days  researchers are also focusing  on the  the optimal additional food to be given to the predators. Some of the findings in this direction show that   the  quality  and quantity of additional food provided  play a crucial role in this optimal studies \cite{srinivasu2018additional}. For instance, the authors in \cite{ananth2021influence,ananth2022achieving} worked on the deterministic Holling type-III predator-prey systems and found the optimal quality and  quantity of additional food required to drive the system to a desired terminal state. \\

It is a well known fact that the real-life systems behave more chaotic.  Stochastic setting can be  an ideal tool to capture such random system dynamics better than that of deterministic system dynamics. Motivated by the above discussions, in this work, we worked on a stochastic predator-prey system with  additional food for predator and exhibiting Holling type-III functional response which in a way can be considered to be a generalization of the work   \cite{ananth2021influence,ananth2022achieving}. We have also incorporated the intra-specific competition among predators in the model to make the model more realistic.  \\

In this work we initially considered a Lagrangian optimal control problem with a cost that is linear w.r.t. state and quadratic w.r.t. control with the end goal of biological conservation of both the species. We considered two cases with the quality of additional food and the quantity of additional food as control variables respectively. We calculated the optimal control values using stochastic maximum principle. We then numerically simulated using Forward-Backward Sampling method. These results showed that biological conservation of the system can be achieved by the optimal control strategy of providing  
 the quality  and the quantity of additional food   high initially and  further reducing them over the time.  \\

Secondly, we studied a time-optimal control problem with the end goal of biological conservation of both the species and also to achieve the goal of pest management in minimum time. In this problem, we worked with a multi-dimensional control involving both the quality and quantity of additional food as control variables. The optimal control values are calculated using the stochastic maximum principle. We also plotted these solutions for two different set of parameters, one applied to biological conservation and the other to pest management. In case of biological conservation, both the controls did not exhibit any switch over the time. In case of pest management, optimal quality of the additional food remained low and constant throughout. The optimal quantity of additional food fluctuated initially and remained high throughout.  \\

The present stochastic optimal control studies can further be improved by incorporating alley effect and even multiple prey and predator species to make the model more realistic. Also, the model presented here does not account for time delay. As diffusion term in this system is independent of control, the results obtained will be similar to that of deterministic case. It will be very interesting to take up problems with control in the diffusion. Since the current systems are of higher orders of non-linearity, it turns out that the numerical methods play a crucial role in understanding the chaotic behaviours. Hence, the study of higher-order numerical methods like Stochastic Runge-Kutta methods will enhance the output.  Also not much work has been done where multiple noise are added simultaneously. We aim to take up these studies in the future. \\

\subsection*{Ethical Approval} 
This research did not require ethical approval.

\subsection*{Funding}
This research was supported by National Board of Higher Mathematics(NBHM), Government of India(GoI) under project grant Time Optimal Control and Bifurcation Analysis of Coupled Nonlinear Dynamical Systems with Applications to Pest Management(02011/11/2021NBHM(R.P)/R$\&$D II/10074).

\subsection*{Conflicts of Interest}
The authors have no conflicts of interest to disclose.

\subsection*{Acknowledgments}
The authors dedicate this paper to the founder chancellor of SSSIHL, Bhagawan Sri Sathya Sai Baba. The corresponding author also dedicates this paper to his loving elder brother D. A. C. Prakash who still lives in his heart.

\printbibliography

\end{document}